\documentclass[a4paper,11pt]{article}
\usepackage{latexsym,amsfonts,amsmath}
\usepackage{latexsym,amsfonts,amsmath,mathrsfs}
\usepackage{graphicx}
\usepackage{color}
\usepackage[mathscr]{eucal}
\usepackage{stmaryrd}
\usepackage{graphicx}
\DeclareMathAlphabet\mathbfcal{OMS}{cmsy}{b}{n}

\usepackage{amsbsy,amsmath}

\def\NN{\mathbb{N}}

\def\RR{\mathbb{R}}

\newcommand{\widebar}[1]{\overline{#1}}

\newcommand{\I}[1]{I_{\{#1\}}}

\def\ds{\displaystyle}

\def\nl{\mbox{} \newline }

\newcounter{hypot}

    \newenvironment{hypot}{\begin{list}
      {\hspace{\labelsep}\bfseries Assumption \Alph{hypot}.}
      {\leftmargin=0pt
       \labelwidth=0cm
       \refstepcounter{hypot}
       \def\makelabel##1{##1}}}{\end{list}}

\newcounter{assump}

\newenvironment{assump}{\begin{list}
      {\hspace{\labelsep}(\Alph{hypot}\arabic{assump})}
      {\leftmargin=1cm
       \labelwidth=1cm
       \usecounter{assump}
       \itshape }}{\end{list}}

\begin{document}

\newtheorem{theorem}{Theorem}[section]
\newtheorem{proposition}[theorem]{Proposition}
\newtheorem{lemma}[theorem]{Lemma}
\newtheorem{corollary}[theorem]{Corollary}
\newtheorem{definition}[theorem]{Definition}
\newtheorem{remark}[theorem]{Remark}
\newtheorem{conjecture}[theorem]{Conjecture}
\newtheorem{assumption}[theorem]{Assumption}
\newtheorem{condition}[theorem]{Condition}
\newcommand{\defi}{\stackrel{\triangle}{=}}

\title{Partially observed optimal stopping problem for discrete-time Markov processes}

\author{ \mbox{ }
\\ \\ B. de Saporta
\\ 
\small Universit\'e de Montpellier, France \\
\small IMAG, CNRS UMR 5149, France \\
\small INRIA Bordeaux Sud Ouest\\
\small e-mail :  benoite.de-saporta@umontpellier.fr
\\
\and
\\ \\ F. Dufour
\thanks{Corresponding author: F. Dufour, INRIA, 200 Avenue de la Vieille Tour, 33405 Talence Cedex, France.}
\\ 
\small Bordeaux INP, France\\
\small IMB, CNRS UMR 5251, France \\
\small INRIA Bordeaux Sud Ouest \\
\small e-mail : dufour@math.u-bordeaux1.fr 
\\
\and
\\ \\ C. Nivot
\\ 
\small INRIA Bordeaux Sud Ouest\\
\small Universit\'e de Bordeaux, France\\
\small IMB, CNRS UMR 5251, France\\
\small e-mail : christophe.nivot@inria.fr
}

\maketitle

\begin{abstract}
This paper is dedicated to the investigation of a new numerical method to approximate the optimal stopping problem for a discrete-time continuous state space Markov chain under partial observations. It is based on a two-step discretization procedure based on optimal quantization. First,we discretize the state space of the unobserved variable by quantizing an underlying reference measure. Then we jointly discretize the resulting approximate filter and the observation process. We obtain a fully computable approximation of the value function with explicit error bounds for its convergence towards the true value fonction.\\


\end{abstract}

%\newpage

%%%%%%%%%%%%%%%%%%%%%%%%%%%%%%%%%%%%%%%%%%%%%%%%%%%%%%%%%%%%%%%%%%%%%%%%%%%%%%%%
%%%%%%%%%%%%%%%%%%%%%%%%%%%%%%%%%%%%%%%%%%%%%%%%%%%%%%%%%%%%%%%%%%%%%%%%%%%%%%%%
\section{Introduction}
\label{Intro}
This paper is dedicated to the investigation of a new numerical method to approximate the optimal stopping problem for a discrete-time continuous state space Markov chain under partial observations. This is known to be a difficult problem, but very important for practical applications. Indeed, the usual approach when dealing with partially observed problems is to introduce the filter or belief process, thus converting the problem into a fully observed one, at the cost of an infinite dimensional state space as the filter process is measure-valued. Thus, there is no straightforward way to discretize the state space of the filter process, and one must often choose a balance between the computational load and the accuracy of the approximation.\\

Unlike the huge literature on discrete state space optimal stopping problems, that on continuous state space is scarce. The most relevant papers addressing this problem are \cite{YZ13,Zhou13,ZFM10}. In \cite{YZ13}, the authors do not propose an approximation of the value function, but only computable upper and lower bounds. They do not require any particular assumptions on the Markov process apart from being simulatable, however they do not provide convergence rates either. Our aim in this paper is more ambitious as we want to construct a numerically tractable approximation of the value function with a bound for the convergence rate.
In \cite{Zhou13}, the authors compute an approximation of the value function based on particle filtering and simulations of the chain trajectories. They assume that all distributions have densities with respect to the Lebesgue measure, that the reward function is convex, but again they don't provide convergence rates for the approximation.
Our approach is more general as it allows the Markov chain kernel have a density with respect to a general product of measures, not necessarily the Lebesgue measure, which may be more relevant for some applications where thresholds are involved, for instance.
In \cite{ZFM10}, the authors propose to parametrize the belief state with the exponential family to dramatically reduce its dimension, but they deal with general control problems for infinite discounted cost and stationary policies that are not suitable for optimal stopping problems.\\

In this paper we propose a new approach, inspired by \cite{pham05} that addresses the optimal stopping problem under partial observation for finite state space chains. The key idea of the authors is to approximate simultaneously the filter and observation processes using a series of quantization grids. 
Optimal quantization is an approximation procedure that replaces a continuous state space variable $X$ by a finite state space one $\widehat{X}$ optimally, in the sense that it minimizes the $L_2$ norm of the difference $|X-\widehat{X}|$, see e.g. \cite{pages04,PPP04} and references therein for more details and applications to numerical probability.
The quantization approach of \cite{pham05} is especially efficient if the state space of the unobserved variable is finite and small. One first simple idea to turn our continuous state space problem into a discrete one is to discretize the state space of the unobserved variables using a regular cartesian grid. However, to ensure precision this may require a huge number of points and possibly useless computations if some areas of the state space are seldom visited.
A better idea is to use the same quantization approach as \cite{pham05} to discretize the unobserved component. This will ensure that the grids have more points in the areas of high density, and are dynamically adapted with time. The state space of the unobserved variable would then be finite, but with time-varying, making the discretization of the filter numerically intractable. 
Our approach attempts at taking the advantages of both these ideas, while minimizing their drawbacks.
Our approximation procedure is in two steps. First, we discretize the state space of the unobserved variable by optimally quantizing an underlying reference distribution. Thus we have a fixed finite state space for the unobserved variable, and the points are optimally distributed to ensure precision at a minimal computational cost.
This yields an approximate filter process that is measure-valued, but can be seen as taking values in a finite dimensional simplex. We then jointly quantize the approximate filter and observation processes. Throughout this procedure, we are able to compute an explicit upper bound for the error, that goes to zero as the number of points in the quantization grids goes to infinity.
\\

The paper is organized as follows. In Section~\ref{sec:pb}, we state the optimal stopping problem under partial observation we are intereted in approximating, and we give the equivalent completely observed sequential decision making problem. In Section~\ref{sec:approx}, we detail our two-step numerical scheme and evaluate the error bound. Section~\ref{sec:example} is dedicated to a numerical example, and the most technical results are postponed to an Appendix.
%%%%%%%%%%%%%%%%%%%%%%%%%%%%%%%%%%%%%%%%%%%%%%%%%%%%%%%%%%%%%%%%%%%%%%%%%%%%%%%%
%%%%%%%%%%%%%%%%%%%%%%%%%%%%%%%%%%%%%%%%%%%%%%%%%%%%%%%%%%%%%%%%%%%%%%%%%%%%%%%%
\section{Problem formulation}\label{sec:pb}
%%%%%%%%%%%%%%%%%%%%%%%
%For \textit{weak} formulation of optimal stochastic control problem, \rose{the reader may consult the references}
%\cite{dufour04,elkaroui87,haussmann90}.
%
We start with some general notation that will be in force throughout the paper. 
$\NN$ is the set of natural numbers including $0$, $\NN^{*}=\NN-\{0\}$, $\RR$ denotes the set of real numbers, $\RR_{+}$ the set of non-negative real numbers,
$\RR_{+}^{*}=\RR_{+}-\{0\}$.
For any $(p,q)\in \NN^{2}$ with $p\leq q$, $\llbracket p ; q \rrbracket$ is the set $\{p,p+1,\ldots,q\}$.
Given $x$ in the Euclidean space $\RR^{n}$,$|x|$ will denote its Euclidean norm. 
Let $I_{E}$ be the indicator function of a set $E$. 
Let $E$ be a metric space where $d$ denotes its associated distance. 
Its Borel $\sigma$-algebra will be denoted by $\mathcal{B}(E)$ and $\mathcal{P}(E)$ is the set of probability measures on $(E,\mathcal{B}(E))$.
The space of all bounded real-valued measurable function on $E$ is denoted by $\mathbb{B}(E)$.
The space $\mathbb{L}(E)$ of all real-valued bounded Lipschitz continuous functions on $E$ is equipped with the norm 
$\|f\|_{\mathbb{L}(E)}=\|f\|_{sup}+L_f$ where $\|f\|_{sup}=\sup_{x\in E} |f(x)|$ and $L_f=\sup_{x\neq y} \frac{|f(x)-f(y)|}{d(x,y)}$ and $\mathbb{L}_{1}(E)=\{f\in \mathbb{L}(E) : \|f\|_{\mathbb{L}(E)} \leq 1\}$.
On $\mathcal{P}(E)$, let us introduce the distance $d_{\mathcal{P}}$ defined by $d_{\mathcal{P}}(\mu,\nu)=\sup_{f\in \mathbb{L}_{1}(E)} \big\{ \int_{E} f d\mu -  \int_{E} f d\nu \big\}$.
The Dirac probability measure concentrated at $x\in E$ will be denoted by $\delta_x$.
If $F$ is a metric space and  $v$ is a real-valued bounded measurable function defined on $E\times F$ and $\gamma$ is a probability measure on $(E,\mathcal{B}(E))$  then by a slight abuse of notation
we write $v(\gamma,y)=\int_{E} v(x,y) \gamma(dx)$ for any $y\in F$.

\subsection{Optimal stopping}
%%%%%%%%%%%%%%%%%%%%%%%
In this section, we describe the optimal stopping problem we are interested in by using a \textit{weak} formulation.
Consider $\mathbb{X} \in \mathcal{B}(\RR^{m})$, $\mathbb{Y} \in \mathcal{B}(\RR^{n})$, a stochastic kernel $R$ on $\mathbb{X}\times \mathbb{Y}$ and a performance function $\mathbf{H}\in \mathbb{B}(\mathbb{X}\times \mathbb{Y})$.

\begin{definition}
\label{Def-opt-stop}
The control is defined by the following term:
\begin{eqnarray*}
\ell & = &\big( \mathbf{\Xi},\mathbfcal{G},\mathbf{Q},\{\mathbfcal{G}_{t}\}_{t\in \llbracket 0 ; N_{0} \rrbracket},\{\mathbfcal{X}_{t},\mathbfcal{Y}_{t}\}_{t\in \llbracket 0 ; N_{0} \rrbracket},\tau\big)
\end{eqnarray*}
\begin{itemize}
\item $\big( \mathbf{\Xi},\mathbfcal{G},\mathbf{Q},\{\mathbfcal{G}_{t}\}_{t\in \llbracket 0 ; N_{0} \rrbracket} \big)$ is a filtered probability space,
\item $\{\mathbfcal{X}_{t},\mathbfcal{Y}_{t}\}_{t\in \llbracket 0 ; N_{0} \rrbracket}$ is an $\mathbb{X}\times\mathbb{Y}$-valued $\{\mathbfcal{G}_{t}\}_{t\in \llbracket 0 ; N_{0} \rrbracket}$-Markov chain defined on $\big( \mathbf{\Xi},\mathbfcal{G},\mathbf{Q} \big)$ where $R$ is its associated transition kernel and $\delta_{(\mathbf{x},\mathbf{y})}$ is its initial distribution,
\item $\tau$ is a $\{\mathbfcal{G}^{\mathbfcal{Y}}_{t}\}_{t\in \llbracket 0 ; N_{0} \rrbracket}$-stopping time where 
$\mathbfcal{G}^{\mathbfcal{Y}}_{t}=\sigma\{\mathbfcal{Y}_{0},\ldots,\mathbfcal{Y}_{t}\}$.
\end{itemize}
\end{definition}
In this setting, $\mathbfcal{X}_{t}$ denotes the hidden variables and $\mathbfcal{Y}_{t}$ the observed variables. Hence the stopping decision $\tau$ depends only on the observations.
The set of the previous controls is denoted by $L$ and $\mathbfcal{E}^{\mathbf{Q}}$ denotes the expectation under the probability $\mathbf{Q}$.
For a control $\ell\in L$, the performance criterion is given by
\begin{eqnarray}
\mathbfcal{H}(\mathbf{x},\mathbf{y},\ell) & =  &\mathbfcal{E}^{\mathbf{Q}}\big[ \mathbf{H}(\mathbfcal{X}_{\tau},\mathbfcal{Y}_{\tau})\big].
\end{eqnarray}
In the previous expression, we write explicitly the dependence of the cost function on the initial state of the Markov chain.
The partially observed optimal stopping problem we are interested in is to maximize the reward function $\mathbfcal{H}(\mathbf{x},\mathbf{y},\ell)$ over $L$. The corresponding value function is thus
\begin{equation}\label{eq:value}
\overline{\mathbfcal{H}}(\mathbf{x},\mathbf{y})=\sup_{\ell\in L} \mathbfcal{H}(\mathbf{x},\mathbf{y},\ell).
\end{equation}
The aim of this paper is to propose a numerical approximation of $\overline{\mathbfcal{H}}(\mathbf{x},\mathbf{y})$ that can be computed in practice and derive bounds for the approximation error.\\

We make the following main assumptions on the parameters of the Markov chain and the performance function. The first ones are mild and state that the transition kernel $R$ of the Markov chain has a density with respect to a product of reference probability measures, and that this density is bounded with Lipschitz regularity. Assumption~\ref{HC} is technical and more restrictive. It states that the density should also be bounded from below. Finally we assume that the performance function is also bounded and Lipschitz-continuous.
\begin{hypot}
\item\label{HA} There exist $\lambda \in \mathcal{P}(\mathbb{X})$, $\nu \in \mathcal{P}(\mathbb{Y})$ and an $\RR_{+}$-valued measurable function $r$ defined on $(\mathbb{X}\times \mathbb{Y})^{2}$
such that
\begin{assump}
\item \label{HA1} for any $(x,y)\in \mathbb{X}\times \mathbb{Y}$, $B\in \mathcal{B}(\mathbb{X})$, $C\in \mathcal{B}(\mathbb{Y})$
$$R(B\times C|x,y)=\int_{B\times C} r(u,v,x,y) \lambda(du)\nu(dv),$$
\item \label{HA2} $\ds \int_{\mathbb{X}} |x|^{2+\beta} \lambda(dx) < \infty$ for some $\beta>0$.
\end{assump}
\end{hypot}

\begin{hypot}
\item\label{HB} There exist positive constants $\widebar{r}$ and $L_{r}$ such that
\begin{assump}
\item \label{HB1} $\ds \sup_{(u,v,x,y)\in (\mathbb{X}\times \mathbb{Y})^{2}} r(u,v,x,y) \leq \widebar{r}$
\item \label{HB2}  for any $(u,v,x,y)\in (\mathbb{X}\times \mathbb{Y})^{2}$, $(u',x',y') \in \mathbb{X}\times \mathbb{X}\times \mathbb{Y}$
$$|r(u,v,x,y)-r(u',v,x',y')|\leq L_{r} \big[ |u-u'|+|x-x'|+|y-y'|\big].$$
\end{assump}
\end{hypot}

\begin{hypot}
\item\label{HC} There exists $\delta>0$ such that $r(\lambda,v,x,y)\geq \delta^{-1}$
for any $(v,x,y) \in \mathbb{Y}\times \mathbb{X}\times \mathbb{Y}$.
\end{hypot}

\begin{hypot}
\item\label{HD} The function $\mathbf{H}$ belongs to $\mathbb{L}(\mathbb{X}\times \mathbb{Y})$.
\end{hypot}
%%%%%%%%%%%%%%%%%%%%%%%%%%%%%%%%%%%%%%%%%%%%%%%%%%%%%%

Our approximation strategy is in three steps. First, we rewrite the problem as a sequential decision-making problem for a fully observed Markov chain on $\mathcal{P}(\mathbb{X})\times\mathbb{Y}$. Then we propose a first approximation based on the discretization of the state space $\mathbb{X}$ by a finite grid $\Gamma_X^N$. Finally, we use a second approximation procedure to discretize the resulting Markov chain on $\mathcal{P}(\Gamma_X^N)\times\mathbb{Y}$.
%%%%%%%%%%%%%%%%%%%%%%
\subsection{Auxiliary completely observed control problem}
\label{aux-compl-obs}
As explained in the introduction, the standard approach to deal with partial observation is to introduce the filter process and convert the problem into a fully observed one on an infinite dimensional state space. In this section, we introduce the auxiliary completely observed control model $\mathcal{M}$. We follow closely the framework of Chapter 5 in \cite{bauerle11}.
The objective of this section is twofold. First, we show that the optimal stopping problem introduced in Definition \ref{Def-opt-stop} is equivalent to a fully observed optimization problem defined in terms of the control model $\mathcal{M}$ below. Second, we prove that the value function $\overline{\mathbfcal{H}}(\mathbf{x},\mathbf{y})$ defined in (\ref{eq:value}) can be obtained by iterating a Bellman operator.\\

As defined in \cite{bauerle11}, let us consider the Bayes' operator $\Phi : \mathbb{Y}\times\mathcal{P}(\mathbb{X})\times\mathbb{Y} \mapsto \mathcal{P}(\mathbb{X})$ given by
$\ds \Phi(v,\theta,y)(du)=\frac{r(u,v,\theta,y)}{r(\lambda,v,\theta,y)}\lambda(du)$
and the stochastic kernel $S$ on $\mathcal{P}(\mathbb{X})\times\mathbb{Y}$ defined by
\begin{eqnarray}
S(B\times C |\theta,y)=\int_{C} \delta_{\Phi(v,\theta,y)}(B) R(\mathbb{X},dv|\theta,y),
\label{Def-S}
\end{eqnarray}
for any $B\in\mathcal{B}(\mathcal{P}(\mathbb{X}))$, $C\in \mathcal{B}(\mathbb{Y})$
and $(\theta,y)\in \mathcal{P}(\mathbb{X})\times\mathbb{Y}$.
For notational convenience, let us introduce the real-valued function $H$ (respectively, $h$) defined on $\mathbb{X}\times\mathbb{Y}\times\{0,1\}\times\{0,1\}$
(respectively, $\mathbb{X}\times\mathbb{Y}\times\{0,1\}$) by
$H(x,y,z,a)=\mathbf{H}(x,y)\I{(z,a)=(0,1)}$ (respectively, $h(x,y,z)=\mathbf{H}(x,y)\I{z=0}$).

Consider the following auxiliary model $\mathcal{M}:= \big(\mathbb{S}, \mathbb{A},Q,H,h\big)$ where
\begin{enumerate}
\item [(a)] the state space is given by $\mathbb{S}=\mathcal{P}(\mathbb{X})\times\mathbb{Y}\times\{0,1\}$,
\item [(b)] the action space is $\mathbb{A}=\{0,1\}$,
\item [(c)] the transition probability function $Q$ is the stochastic kernel on $\mathbb{S}$ given $\mathbb{S}\times\mathbb{A}$ defined by
$Q(B\times C \times D|\theta,y,z,a)=S(B\times C |\theta,y) \big[ \delta_{z}(D) \I{a=0}+ \delta_{1}(D) \I{a=1} \big]$
for any $B\in\mathcal{B}(\mathcal{P}(\mathbb{X}))$, $C\in \mathcal{B}(\mathbb{Y})$, $D\subset\{0,1\}$ and $(\theta,y,z,a)\in \mathbb{S}\times\mathbb{A}$,
\item [(d)] the cost-per-stage is $H(\theta,y,z,a)$ and the terminal cost is $h(\theta,y,z)$
for any $(\theta,y,z,a)\in \mathbb{S}\times\mathbb{A}$ (recalling the slight abuse of notation introduced at the end of Section \ref{Intro}).
\end{enumerate}

The underlying idea is that the filtered trajectory is constructed recursively thanks to the Bayes operator $\Phi$, and the kernel $S$ is the driving kernel of the Markov chain of the filter and observations. The optimal stopping problem is then stated as a sequential decision making problem where at each time step the controller may stop (action $a=1$) or continue (action $a=0$). The additional variable $z\in\{0,1\}$ indicates whether the trajectory has already been stopped ($z=1$) or not ($z=0$).

Introduce $\Omega=\mathbb{S}^{N_{0}+1}$, $\mathcal{F}$ its associated product $\sigma$-algebra and the coordinate projections
$\Theta_{t}$ (respectively $Y_{t}$ and $Z_{t}$) from $\Omega$ to the set $\mathcal{P}(\mathbb{X})$ (respectively $\mathbb{Y}$ and $\{0,1\}$).
Let $\Pi^{o}$ be the set of all deterministic past dependent control policies $\pi=\{\pi_{t}\}_{t\in \llbracket 0 ; N_{0}-1 \rrbracket}$ where
$\pi_{0}$ is a measurable $\mathbb{A}$-valued function defined on $\mathbb{Y}\times\{0,1\}$ and 
$\pi_{t}$ is a measurable $\mathbb{A}$-valued function defined on $(\mathbb{Y}\times\{0,1\}\times \mathbb{A})^{t} \times \mathbb{Y}\times\{0,1\} $ for
$t\in \llbracket 1 ; N_{0} \rrbracket$. 

Consider an arbitrary policy $\pi\in\Pi^{o}$. Define the action process $\{A_{t}\}_{t\in \llbracket 0 ; N_{0}-1 \rrbracket}$ by
$A_{t}=\pi_{t}(Y_{0},Z_{0},A_{0},\ldots,Y_{t-1},Z_{t-1},A_{t-1},Y_{t},Z_{t})$
for $t\in \llbracket 1 ; N_{0}-1 \rrbracket$ and $A_{0}=\pi_{0}(Y_{0},Z_{0})$.
Define $\mathcal{F}_{t}=\sigma\{\Theta_{0},Y_{0},Z_{0},\ldots,\Theta_{t},Y_{t},Z_{t}\}$
for $t\in \llbracket 0 ; N_{0} \rrbracket$.
According to \cite{bauerle11,hernandez96}, for an arbitrary policy $\pi\in\Pi^{o}$ there exists a probability measure $P^{\pi}_{(\mathbf{x},\mathbf{y})}$
on $\big(\Omega,\mathcal{F} \big)$ which satisfy
\begin{itemize}
\item[i)] $P^{\pi}_{(\mathbf{x},\mathbf{y})}\big((\Theta_{0},Y_{0},Z_{0})\in B\times C \big)=\delta_{\delta_{\mathbf{x}}}(B) \delta_{\mathbf{y}}(C) \delta_{0}(D)$,
\item[ii)] $P^{\pi}_{(\mathbf{x},\mathbf{y})}\big((\Theta_{t+1},Y_{t+1},Z_{t+1})\in B\times C\times D |\mathcal{F}_{t}\big)
=Q(B\times C\times D|\Theta_{t},Y_{t},Z_{t},A_{t})$,
\end{itemize}
for any $B\in \mathcal{B}(\mathcal{P}(\mathbb{X}))$, $C\in \mathcal{B}(\mathbb{Y})$, $D\subset\{0,1\}$, and $t\in  \llbracket 0 ; N_{0}-1 \rrbracket$.

\bigskip

\noindent
The expectation under the probability $P^{\pi}_{(\mathbf{x},\mathbf{y})}$ is denoted by $E^{\pi}_{(\mathbf{x},\mathbf{y})}$.
For a policy $\pi\in \Pi^{o}$, the performance criterion is given by
\begin{eqnarray}
\mathbfcal{H}_{\mathcal{M}}(\mathbf{x},\mathbf{y},\pi) & =  &
E^{\pi}_{(\mathbf{x},\mathbf{y})} \Big[ \sum_{t=0}^{N_{0}-1}  H(\Theta_{t},Y_{t},Z_{t},A_{t}) \Big]
+E^{\pi}_{(\mathbf{x},\mathbf{y})}\big[ h(\Theta_{N_{0}},Y_{N_{0}},Z_{N_{0}}) \big].
\end{eqnarray}
The optimization problem consists in maximizing the reward function $\mathbfcal{H}_{\mathcal{M}}(\mathbf{x},\mathbf{y},\pi)$ over $\Pi^{o}$ and the corresponding value function is
\begin{equation}\label{eq:valueM}
\overline{\mathbfcal{H}}_{\mathcal{M}}(\mathbf{x},\mathbf{y})=\sup_{\pi\in \Pi^{o}} \mathbfcal{H}_{\mathcal{M}}(\mathbf{x},\mathbf{y},\pi).
\end{equation}
It can be computed using dynamic programming. Consider the Bellman operator $\mathfrak{B}$ defined on $\mathbb{B}( \mathcal{P}(\mathbb{X})\times\mathbb{Y})$ by
\begin{equation}\label{def:Bellman}
\mathfrak{B}f(\theta,y)=\max\{\mathbf{H}(\theta,y),Sf(\theta,y)\},
\end{equation}
for $f\in \mathbb{B}( \mathcal{P}(\mathbb{X})\times\mathbb{Y})$. It should be clear that under Assumption \ref{HD},  $\mathfrak{B}$ maps $\mathbb{L}( \mathcal{P}(\mathbb{X})\times\mathbb{Y})$ onto 
$\mathbb{B}( \mathcal{P}(\mathbb{X})\times\mathbb{Y})$.
For notational convenience, $\mathfrak{B}^{k} $ denotes the $k$-th iteration of $\mathfrak{B}$ recursively defined by $\mathfrak{B}^{0}f=f$, $\mathfrak{B}^{1}f=\mathfrak{B}f$ and $\mathfrak{B}^{k}f=\mathfrak{B}(\mathfrak{B}^{k-1}f)$
for $k\in \llbracket 2 ; N_{0}\rrbracket$ and $f\in \mathbb{B}(\mathcal{P}(\mathbb{X})\times\mathbb{Y})$.
\begin{theorem}
\label{Equivalence-Costs}
Suppose Assumptions \ref{HA1}, \ref{HB} and \ref{HD} hold. Then
%The following assertions hold.
%\nl
%$i)$ For any control $\ell \in L$, there exist a policy $\pi\in \Pi^{o}$ such that
%$$\mathbfcal{H}_{M}(\mathbf{x},\mathbf{y},\pi)=\mathbfcal{H}(\mathbf{x},\mathbf{y},\lambda).$$
%$ii)$ For any policy $\pi\in \Pi^{o}$, there exist a control $\ell \in L$ such that
%$$\mathbfcal{H}(\mathbf{x},\mathbf{y},\lambda)=\mathbfcal{H}_{M}(\mathbf{x},\mathbf{y},\pi).$$
%Therefore,
\begin{eqnarray}
\overline{\mathbfcal{H}}(\mathbf{x},\mathbf{y})=\overline{\mathbfcal{H}}_{\mathcal{M}}(\mathbf{x},\mathbf{y})=\mathfrak{B}^{N_{0}}\mathbf{H}(\delta_{\mathbf{x}},\mathbf{y}).
\label{Cost=}
\end{eqnarray}
%Define recursively the sequence $\{V_{k}\}_{ k \in \llbracket 0 ; N_{0}\rrbracket}$ of real-valued bounded measurable functions on $\mathcal{P}(\mathbb{X})\times\mathbb{Y}$ by
%$V_{0}(\theta,y)=\mathbf{H}(\theta,y)$ for $(\theta,y)\in \mathcal{P}(\mathbb{X})\times\mathbb{Y}$ and
%$V_{k}=\mathfrak{B}V_{k-1}$ for $k \in \llbracket 1 ; N_{0}\rrbracket$ where the Bellman operator $\mathfrak{B}$ is defined by
%$$\mathfrak{B}f(\theta,y)=\max\{\mathbf{H}(\theta,y),Sf(\theta,y)\}.$$
%Then $\overline{\mathbfcal{H}}(\mathbf{x},\mathbf{y})=V_{N_{0}}(\delta_{\mathbf{x}},\mathbf{y})$.
\end{theorem}
\textbf{Proof:}
See Appendix \ref{aux-part-obs}. \hfill $\Box$

%%%%%%%%%%%%%%%%%%%%%%%%%%%%%%%
\section{Approximation results}\label{sec:approx}
%%%%%%%%%%%%%%%%%%%%%%%%%%%%%%%
We now build our approximation procedure for the value function $\overline{\mathbfcal{H}}_{\mathcal{M}}(\mathbf{x},\mathbf{y})$. It is based on two discretizations. First, we replace the continuous state space $\mathbb{X}$ of the hidden variable by a discrete one $\Gamma^N_X$ with cardinal $N$. 
Thus, model $\mathcal{M}$ can be approximated by a similar sequential decision making problem for a Markov chain on the finite dimensional state space $\mathcal{P}(\Gamma^{N}_{X})\times\mathbb{Y}$. We then discretize the latter Markov chain using time-dependent grids following the same procedure as in \cite{pham05}.
In both steps, the discretization grids we use are quantization grids. They are especially appealing because they are optimized so that there are more points in the areas of high density, and they allow to control the discretization error in $L_2$-norm as long as the underlying operators have Lipschitz continuity properties.
In this section, we first recall the basics of optimal quantization, then present the two discretization steps and derive the discretization error.

%%%%%%%%%%%%%%%%%%%%%%%%%%%%%%%
\subsection{Optimal quantization}\label{sec:quantize}
%%%%%%%%%%%%%%%%%%%%%%%%%%%%%%%

Consider an $\mathbb{R}^{d}$-valued random variable $Z$ defined on a probability space $(G,\mathcal{G},\mathbb{P})$ (with corresponding expectation operator
$\mathbb{E}$) such that $\|Z\|_2 < \infty$ where $\| Z \|_2$ denotes the $L_{2}$-norm of $Z$.
Let $N$ be a fixed integer. The optimal $L_{2}$-quantization of the random variable $Z$ consists in finding the best possible $L_{2}$-approximation of $Z$
by a random variable $\widehat{Z}_{N}$ on $(G,\mathcal{G},\mathbb{P})$ taking at most $N$ values in $\RR^{d}$,
which will be denoted by $\{z^{1}_{N},\ldots,z^{N}_{N}\}$.
The asymptotic properties of the $L_{2}$-quantization are summarized in the following result (see, e.g., \cite[Theorem 3]{bally05}), which uses the notation $p_{\Gamma}(z)$ for the closest neighbor projection of $z\in\RR^{d}$ on a grid $\Gamma=\{z_{1},\ldots,z_{N}\}\subseteq\RR^{d}$.
\begin{theorem}\label{th-quantization}
Let $Z$ be an $\RR^{d}$-valued random variable on $(G,\mathcal{G},\mathbb{P})$, and suppose that for some $\epsilon>0$ we have $\mathbb{E}[|Z|^{2+\epsilon}]<+\infty$. Then
\begin{eqnarray*}
\lim_{N\rightarrow \infty} N^{2/d} \min_{|\Gamma|\leq N} \| Z-p_{\Gamma}(Z)\|^{2}_{2}& =& J_{d,2} \int_{\RR^{d}} |h_{Z}(u)|^{d/(2+d)}(u) du,
\end{eqnarray*}
where $h_{Z}(u)$ denotes the density of the absolutely continuous part of the distribution of $Z$ with respect to the Lebesgue measure on $\mathbb{R}^{d}$, and $J_{2,d}$ is a universal constant.
\end{theorem}
Finally, let us mention that there exist algorithms that can numerically find, for a fixed $N$, the quantization of~$Z$ (or, equivalently, the grid $\{z^{1}_{N},\ldots,z^{N}_{N}\}$ attaining the minimum in Theorem \ref{th-quantization} above and its distribution) as soon as $Z$ is simulatable. Basically, the quantization grids for the variable $Z$ will have more points in the areas of high density, and fewer points in the areas of low density for $Z$.

%%%%%%%%%%%%%%%%%%%%%%%%%%%%%%%
\subsection{First approximation}
%%%%%%%%%%%%%%%%%%%%%%%%%%%%%%%
\label{approx-1}
The main originality of our work is to propose a discretization of the state space $\mathbb{X}$ based on the quantization of the reference measure $\lambda$ defined in Assumption~\ref{HA1}. It greatly helps minimizing the computational burden as only one grid is required, instead of a series of grids as one usually does when trying to quantize accurately a Markov chain. In addition, we obtain bounds for the error. However this come at a cost: in order to guarantee that the approximated transition kernel is still a Markov kernel, the density $r$ also appears in the denominator. This is why we need the lower bound of Assumption~\ref{HC} to control the error.\\

To build the approximation and evaluate the error thus entailed, we first quantize the reference probability $\lambda$. Then we replace it by its quantized approximation $\lambda_N$ in the definitions of kernels $R$, $S$, the Bayes operator $\Phi$ and plug these approximations into the Bellman operator. We obtain an approximate Bellman operator $\mathfrak{B}_N$ and our first approximation of the value function is build by iterating $\mathfrak{B}_N$, following Equation (\ref{Cost=}).\\

According to the previous discussion of Section~\ref{sec:quantize}, given an integer $N$, let $\widehat{X}_{N}$ be the optimal $L_{2}$-quantization of the random variable $X$ with distribution $\lambda$
on a probability space $(G,\mathcal{G},\mathbb{P})$ ($\mathbb{E}[\cdot]$ will stand for the expectation associated to $\mathbb{P}$).
Let us denote by $\Gamma^{N}_{X}=\{x^{1}_{N},\ldots,x^{N}_{N}\}$ an optimal grid. There is no loss of generality to assume that $x^{1}_{N}=\mathbf{x}$.
We write $\lambda_{N}$ for the distribution of $\widehat{X}_{N}$, that is, $\lambda_{N}(du)=\sum_{i=1}^{N}\mathbb{P}(\widehat{X}_{N}=x^{i}_{N})\delta_{x_N^i}(du)$
and 
$$\epsilon_{N}=\|X-\widehat{X}_{N}\|_{2}$$ 
for the $L_{2}$-quantization error between $X$ and $\widehat{X}_{N}$. 
Assume also the existence of a random variable $Y$ with distribution $\nu$ on $(G,\mathcal{G},\mathbb{P})$.\\

We define the quantized approximations of kernels $R$ and $S$ by plugging-in $\lambda_N$ as follows.
Consider the stochastic kernel $R_{N}$ on $\mathbb{X}\times\mathbb{Y}$ where
\begin{eqnarray}
R_{N}(B\times C|x,y)=\int_{B\times C} \frac{r(u,v,x,y)}{r(\lambda_{N},\nu,x,y)} \lambda_{N}(du)\nu(dv),
\label{Def-R-N}
\end{eqnarray}
for any $(x,y)\in \mathbb{X}\times \mathbb{Y}$, $B\in \mathcal{B}(\mathbb{X})$, $C\in \mathcal{B}(\mathbb{Y})$.
Note that the support of $R_{N}(\cdot|x,y)$ is actually $\Gamma^{N}_{X}\times \mathbb{Y}$.
Let us introduce the stochastic kernel $S_{N}$ on $\mathcal{P}(\mathbb{X})\times\mathbb{Y}$ defined by
\begin{eqnarray}
S_{N}(B\times C |\theta,y)=\int_{C} \delta_{\Phi_{N}(v,\theta,y)}(B) R_{N}(\mathbb{X},dv|\theta,y),
\label{Def-S-N}
\end{eqnarray}
where $\Phi_{N} : \mathbb{Y}\times\mathcal{P}(\mathbb{X})\times\mathbb{Y} \to \mathcal{P}(\mathbb{X})$ given by
\begin{eqnarray}
\Phi_{N}(v,\theta,y)(du)=\frac{r(u,v,\theta,y)}{r(\lambda_{N},v,\theta,y)}\lambda_{N}(du).
\label{Def-Phi-N}
\end{eqnarray}
Here again, $\Phi_{N}$ actually maps $\mathbb{Y}\times\mathcal{P}(\mathbb{X})\times\mathbb{Y}$ onto $\mathcal{P}(\Gamma^{N}_{X})$ and the support of $S_{N}(\cdot  |\theta,y)$ is $\mathcal{P}(\Gamma^{N}_{X})\times \mathbb{Y}$.
Next we define the approximated Bellman operator. Consider the operator $\mathfrak{B}_{N}$ defined on $\mathbb{B}( \mathcal{P}(\mathbb{X})\times\mathbb{Y})$ by
\begin{eqnarray}
\mathfrak{B}_{N}f(\theta,y)=\max\{\mathbf{H}(\theta,y),S_{N}f(\theta,y)\}
\label{Def-B-N}
\end{eqnarray}
for $f\in \mathbb{B}( \mathcal{P}(\mathbb{X})\times\mathbb{Y})$. It should be clear that under Assumption \ref{HD},  $\mathfrak{B}_{N}$ maps $\mathbb{B}( \mathcal{P}(\mathbb{X})\times\mathbb{Y})$ onto 
$\mathbb{B}( \mathcal{P}(\mathbb{X})\times\mathbb{Y})$ and $\mathbb{B}( \mathcal{P}(\Gamma^{N}_{X})\times\mathbb{Y})$ onto 
$\mathbb{B}( \mathcal{P}(\Gamma^{N}_{X})\times\mathbb{Y})$.
For notational convenience, $\mathfrak{B}_{N}^{k} $ denotes the $k$-th iteration of $\mathfrak{B}_{N}$.
The first approximate value function is then defined as the $N_0$-iterate of $\mathbb{B}_N$ on $\mathbf{H}$:
\begin{equation}\label{value:N}
\overline{\mathbfcal{H}}_{\mathcal{M},N}(\mathbf{x},\mathbf{y})=\mathfrak{B}_{N}^{N_{0}}\mathbf{H}(\delta_{\mathbf{x}},\mathbf{y}).
\end{equation}
We now study the error induced by this discretization on the value function. To do so, we study how the quantization error is propagated through the Bellman operator. We start with technical results on the density of the Markov chain integrated with respect to $\lambda$ and with respect to its quantized approximation $\lambda_N$.
\begin{lemma}
\label{bound-1/r-N}
For any $(\theta,y)\in \mathcal{P}(\mathbb{X})\times\mathbb{Y}$ and $N\in \NN$ such that
$\displaystyle \epsilon_{N}\leq \frac{1}{2L_{r}}$, one has
\begin{align*}
\frac{1}{r(\lambda_{N},\nu,\theta,y)}\leq 2 \quad \mbox{ and } \quad \Big| 1-\frac{1}{r(\lambda_{N},\nu,\theta,y)} \Big| \leq 2 L_{r} \epsilon_{N}.
\end{align*}
\end{lemma}
\textbf{Proof:}
Clearly, we have $\big| 1 - r(\lambda_{N},\nu,\theta,y) \big| = \big| r(\lambda,\nu,\theta,y) - r(\lambda_{N},\nu,\theta,y) \big| \leq  L_{r} \epsilon_{N}$ by using Assumption
\ref{HB2} and so, $\displaystyle \frac{1}{r(\lambda_{N},\nu,\theta,y)}\leq 2$ for $N\in \NN$ satisfying $\displaystyle \epsilon_{N}\leq \frac{1}{2L_{r}}$. Therefore,
\begin{align*}
\Big| 1-\frac{1}{r(\lambda_{N},\nu,\theta,y)} \Big| & \leq \frac{1}{r(\lambda_{N},\nu,\theta,y)} \big| r(\lambda,\nu,\theta,y) - r(\lambda_{N},\nu,\theta,y) \big|
\leq \frac{1}{r(\lambda_{N},\nu,\theta,y)} L_{r} \epsilon_{N},
\end{align*}
giving the result.
\hfill $\Box$

\begin{lemma}
\label{bound-1/r-N-bis}
Suppose Assumption \ref{HC} holds. For any $(\theta,y)\in \mathcal{P}(\mathbb{X})\times\mathbb{Y}$ and $N\in \NN$ such that
$\displaystyle \epsilon_{N}\leq \frac{1}{2\delta L_{r}}$, one has
$\ds \frac{1}{r(\lambda_{N},v,\theta,y)}\leq 2\delta$.
\end{lemma}
\textbf{Proof:}
The proof is similar to the one of Lemma \ref{bound-1/r-N} and is therefore, omitted.
\hfill $\Box$

\begin{lemma}
\label{inv-r-inv-r-N}
Suppose Assumptions \ref{HB} and \ref{HC} hold. For any $\theta\in \mathcal{P}(\mathbb{X})$ and $N\in \NN$ such that $\epsilon_{N}\leq \frac{1}{2\delta L_{r}}$, one has
\begin{align*}
\sup_{y\in \mathbb{Y}} d_{\mathcal{P}}\big(\Phi(Y,\theta,y),\Phi_{N}(Y,\theta,y)\big) \leq
\delta \big[(2\delta+1) L_{r} +\widebar{r} \big] \epsilon_{N}.
\end{align*}
\end{lemma}
\textbf{Proof:} Consider $f\in \mathbb{L}_{1}(\mathbb{X})$. We have
\begin{align*}
\Big| \int_{\mathbb{X}} f(x) \Phi(v,\theta,y)(dx) - \int_{\mathbb{X}} & f(x) \Phi_{N}(v,\theta,y)(dx) \Big|
\leq \mathbb{E} \Big[ \Big| f(X) \frac{r(X,v,\theta,y)}{r(\lambda,v,\theta,y)} -f(X_{N}) \frac{r(X_{N},v,\theta,y)}{r(\lambda_{N},v,\theta,y)} \Big| \Big]
\nonumber \\
& \leq \mathbb{E} \Big[ |f(X_{N})| r(X_{N},v,\theta,y) \Big| \frac{1}{r(\lambda,v,\theta,y)} - \frac{1}{r(\lambda_{N},v,\theta,y)} \Big| \Big] 
\nonumber \\
& \phantom{\leq} + \mathbb{E} \Big[  \frac{|f(X)|}{r(\lambda,v,\theta,y)} \big| r(X,v,\theta,y) - r(X_{N},v,\theta,y) \big| \Big]
\nonumber \\
& \phantom{\leq} + \mathbb{E} \Big[ \frac{r(X_{N},v,\theta,y)}{r(\lambda,v,\theta,y)}  | f(X)  -f(X_{N}) | \Big].
\end{align*}
By using Assumptions \ref{HB1} and \ref{HC}, it follows that
\begin{align}
\Big| \int_{\mathbb{X}} f(x) & \Phi(v,\theta,y)(dx) - \int_{\mathbb{X}} f(x) \Phi_{N}(v,\theta,y)(dx) \Big|
\nonumber \\
& \leq \widebar{r} \mathbb{E} \Big[ \Big| \frac{1}{r(\lambda,v,\theta,y)} -  \frac{1}{r(\lambda_{N},v,\theta,y)} \Big| \Big] 
+ \delta L_{r}\epsilon_{N} +\delta \widebar{r} \epsilon_{N}.
\label{eq1-inv-r-inv-r-N}
\end{align}
However, for $N\in \NN$ satisfying  $\epsilon_{N}\leq \frac{1}{2\delta L_{r}}$ we get from Lemma \ref{bound-1/r-N-bis} that
\begin{align*}
\Big| \frac{1}{r(\lambda,v,\theta,y)} -  \frac{1}{r(\lambda_{N},v,\theta,y)} \Big| \leq \frac{1}{r(\lambda,v,\theta,y)r(\lambda_{N},v,\theta,y)} L_{r} \epsilon_{N}
\leq 2\delta^{2}L_{r}\epsilon_{N}
\end{align*}
and with equation (\ref{eq1-inv-r-inv-r-N}), this shows the result.
\hfill $\Box$

We now need to ensure that both $\mathfrak{B}$ and $\mathfrak{B}_N$ operate on $\mathbb{L}(\mathcal{P}(\mathbb{X})\times\mathbb{Y})$.

\begin{lemma}
\label{B Lip}
Suppose Assumptions \ref{HB}, \ref{HC} and \ref{HD} hold. For any $\theta,\theta'\in \mathcal{P}(\mathbb{X})$, $y,y'\in \mathbf{Y}$, $N\in \NN$ and $f\in \mathbb{L}(\mathcal{P}(\mathbb{X})\times\mathbb{Y})$ one has $\mathfrak{B}f\in \mathbb{L}(\mathcal{P}(\mathbb{X})\times\mathbb{Y})$.
% with
%\begin{align*}
%\|\mathfrak{B}f\|_{ \mathbb{L}(\mathcal{P}(\mathbb{X})\times\mathbb{Y})}&\leq ???\\
%\|\mathfrak{B}_Nf\|_{ \mathbb{L}(\mathcal{P}(\mathbb{X})\times\mathbb{Y})}&\leq ???\\
%\end{align*}
\end{lemma}
\textbf{Proof:} Consider $f \in \mathbb{L}(\mathcal{P}(\mathbb{X}) \times \mathbb{Y})$ and $(\theta,y),(\theta',y') \in \mathcal{P}(\mathbb{X}) \times \mathbb{Y}$. On the one hand,
\begin{align}
|\mathfrak{B} f(\theta,y)| &\leq \max\{|\mathbf{H}(\theta,y)| ; |S f (\theta,y)|\} \nonumber\\
&\leq \max\left\{\|\mathbf{H}\|_{\mathbb{L}(\mathcal{P}(\mathbb{X}) \times \mathbb{Y})} ; \displaystyle\int_{\mathbb{X}\times\mathbb{Y}} \left|f\big(\Phi(y',\theta,y),y'\big)\right|r(x',y',\theta,y)\lambda(dx')\nu(dy')\right\}\nonumber\\
&\leq \max\left\{\|\mathbf{H}\|_{\mathbb{L}(\mathcal{P}(\mathbb{X}) \times \mathbb{Y})} ; \|f\|_{\mathbb{L}(\mathcal{P}(\mathbb{X}) \times \mathbb{Y})}\right\}.
\end{align}

On the other hand,
\begin{align}
|S f(\theta,y) - S f(\theta',y')| &\leq \mathbb{E}\left[\left|f\big(\Phi(Y,\theta,y), Y \big)\right| \left|r(X,Y,\theta,y) - r(X,Y,\theta',y')\right| \right] \nonumber\\
&~~~ + \mathbb{E}\left[ r(X,Y,\theta',y') \left|f\big(\Phi(Y,\theta,y),Y \big) - f\big(\Phi(Y,\theta',y'),Y \big)\right| \right]\nonumber\\
&\leq \|f\|_{\mathbb{L}(\mathcal{P}(\mathbb{X})\times\mathbb{Y})}(\overline{r} + L_r)[d_{\mathcal{P}}(\theta,\theta') + |y - y'|]\nonumber\\
&~~~ + \overline{r}\|f\|_{\mathbb{L}(\mathcal{P}(\mathbb{X}) \times \mathbb{Y})}\mathbb{E}\left[d_{\mathcal{P}}(\Phi(Y,\theta,y),\Phi(Y,\theta',y'))\right].
\end{align}
Let now $g \in \mathbb{L}_1(\mathcal{P}(\mathbb{X}) \times \mathbb{Y})$ and $v \in \mathbb{Y}$. Then, one has
\begin{align}
& \left|\displaystyle \int_{\mathbb{X}} g(u)\Phi(v,\theta,y)(du) - \int_{\mathbb{X}} g(u)\Phi(v,\theta',y')(dy)\right|\nonumber\\
&\leq \mathbb{E}\left[|g(X)|\left|\dfrac{r(X,v,\theta,y)}{r(\lambda,v,\theta, y)} - \dfrac{r(X,v,\theta',y')}{r(\lambda,v,\theta',y')}\right| \right] \nonumber\\
&\leq \nonumber \mathbb{E}\left[\dfrac{1}{r(\lambda,v, \theta, y)} |r(X,v, \theta,y) - r(X,v,\theta',y')| \right]\\
&~~~ + \mathbb{E}\left[ r(X,v,\theta',y') \left|\dfrac{1}{r(\lambda,v,\theta, y)} - \dfrac{1}{r(\lambda,v,\theta',y')}\right| \right]\nonumber\\
&\leq \delta(\overline{r} + L_r)(1 +\delta\overline{r}) [d_{\mathcal{P}}(\theta,\theta') + |y-y'|]
\end{align}
by using assumptions \ref{HB} and \ref{HC}. So, one has
\begin{equation}
\mathbb{E}\left[d_{\mathcal{P}}(\Phi(Y,\theta,y),\Phi(Y,\theta',y'))\right] \leq \delta(\overline{r} + L_r)(1 +\delta\overline{r}) [d_{\mathcal{P}}(\theta,\theta') + |y-y'|].
\end{equation}
Then, by using assumption \ref{HD}, it is straightforward to write
\begin{align}
|\mathfrak{B} f(\theta,y) - \mathfrak{B}f(\theta',y')| &\leq |\mathbf{H}(\theta,y) - \mathbf{H}(\theta',y')| + |S f(\theta,y) - S f(\theta',y')| \nonumber\\
&\leq L_{\mathfrak{B}f}[d_{\mathcal{P}}(\theta,\theta') + |y-y'|]
\end{align}
with
\begin{equation}
L_{\mathfrak{B}f} = \|\mathbf{H}\|_{\mathbb{L}(\mathcal{P}(\mathbb{X}) \times \mathbb{Y})} + \|f\|_{\mathbb{L}(\mathcal{P}(\mathbb{X}) \times \mathbb{Y})}(\overline{r} + L_r)\big(1 + \delta\overline{r}(1 + \overline{r}\delta)\big).
\end{equation}
Thus, $\mathfrak{B}f$ is bounded and Lipschitz-continuous on $\mathcal{P}(\mathbb{X}) \times \mathbb{Y}$.

%Le lemme \ref{Philip} permet de majorer $\E\left[\norme{\Phi(Y,\rho,y) - \Phi(Y,\rho',y')}_{KR}\right]$ et donc d'avoir
%\begin{equation}
%\vabs{S f(\rho,y) - S f(\rho',y')} \leq \norme{f}_{\mbb{L}}(\ov{q} + L_q)\big(1 + \ov{q}\delta(1 + \ov{q}\delta)\big)[\norme{\rho - \rho'}_{KR} + \vabs{y - y'}].
%\end{equation}
%Donc en utilisant la proposition \ref{fctR} et l'inégalité (\ref{ineg}),
%\begin{align}
%\vabs{\B f(\rho,y) - \B f(\rho',y')} &\leq \vabs{R(\rho,y) - R(\rho',y')} + \vabs{S f(\rho,y) - S f(\rho',y')} \nonumber\\
%&\leq \Big(\norme{R}_{\mbb{L}} + \norme{f}_{\mbb{L}}(\ov{q} + L_q)\big(1 + \delta\ov{q}(1 + \ov{q}\delta)\big)\Big)[\norme{\rho - \rho'}_{KR} + \vabs{y - y'}]
%\end{align}
%d'où le résultat.
\hfill $\Box$

\begin{proposition}
\label{Prop-Diff-B-B-N}
Suppose Assumptions \ref{HA}, \ref{HB} and \ref{HC} hold. Let $N\in \NN$ satisfying
$\epsilon_{N}\leq \frac{1}{2L_{r}}(1\wedge \frac{1}{\delta})$,
\begin{eqnarray}
\big| \mathfrak{B}f(\theta,y) - \mathfrak{B}_{N}f(\theta,y)\big| \leq \|f\|_{\mathbb{L}(\mathcal{P}(\mathbb{X})\times\mathbb{Y})}  K_{1} \epsilon_{N}
\label{Diff-B-B-N}
\end{eqnarray}
for any $(\theta,y)\in \mathcal{P}(\mathbb{X})\times\mathbb{Y}$ and $f\in \mathbb{L}(\mathcal{P}(\mathbb{X})\times\mathbb{Y})$ with
\begin{eqnarray}
K_{1}= L_{r}(1+2\widebar{r}L_{r}) +\delta \widebar{r} \big[\widebar{r}+L_{r}(2+2\delta \widebar{r}) \big].
\end{eqnarray}
\end{proposition}
\textbf{Proof:} Consider $(\theta,y)\in \mathcal{P}(\mathbb{X})\times\mathbb{Y}$ and $f\in \mathbb{L}(\mathcal{P}(\mathbb{X})\times\mathbb{Y})$.
Clearly, we have $$\big| \mathfrak{B}f(\theta,y) - \mathfrak{B}_{N}f(\theta,y)\big| \leq \big| Sf(\theta,y) - S_{N}f(\theta,y)\big|.$$
By using the definition of $S$ and $S_{N}$ (see equations (\ref{Def-S}) and (\ref{Def-S-N}) respectively), we have
\begin{align*}
\big| Sf(\theta,y) - S_{N}f(\theta,y)\big| 
& \leq  \mathbb{E} \bigg[\Big| f(\Phi(Y,\theta,y),Y) r(\lambda,Y,\theta,y)- f(\Phi_{N}(Y,\theta,y),Y) \frac{r(\lambda_{N},Y,\theta,y)}{r(\lambda_{N},\nu,\theta,y)} \Big|\bigg]
\nonumber \\
& \leq  \mathbb{E} \bigg[ r(\lambda,Y,\theta,y) \Big| f(\Phi(Y,\theta,y),Y) - f(\Phi_{N}(Y,\theta,y),Y) \Big|\bigg] \nonumber \\
& \phantom{\leq} + \mathbb{E} \bigg[  \big| f(\Phi_{N}(Y,\theta,y),Y) \big| \big| r(\lambda,Y,\theta,y) - r(\lambda_{N},Y,\theta,y) \big| \bigg] \nonumber \\
& \phantom{\leq} + \mathbb{E} \bigg[ r(\lambda_{N},Y,\theta,y) \big| f(\Phi_{N}(Y,\theta,y),Y) \big|  \Big| 1-\frac{1}{r(\lambda_{N},\nu,\theta,y)} \Big| \bigg] .
\end{align*}
Consequently, it follows that
\begin{align*}
\big| \mathfrak{B}f(\theta,y) - \mathfrak{B}_{N}f(\theta,y) \big| 
& \leq  \widebar{r} \|f\|_{\mathbb{L}(\mathcal{P}(\mathbb{X})\times\mathbb{Y})} \mathbb{E} \Big[ d_{\mathcal{P}}\big(\Phi(Y,\theta,y),\Phi_{N}(Y,\theta,y)\big) \Big] \nonumber \\
& \phantom{\leq} + \|f\|_{\mathbb{L}(\mathcal{P}(\mathbb{X})\times\mathbb{Y})} L_{r} \epsilon_{N}+ \widebar{r} \|f\|_{\mathbb{L}(\mathcal{P}(\mathbb{X})\times\mathbb{Y})}
\Big| 1-\frac{1}{r(\lambda_{N},\nu,\theta,y)} \Big|.
\end{align*}
By using Lemma \ref{inv-r-inv-r-N} and \ref{bound-1/r-N}, we get the result.
\hfill $\Box$\\

We can now state and prove the main result of this section bounding the error between the true value function and its quantized approximation.
\begin{theorem}
\label{Theo-approx-1}
Suppose Assumptions \ref{HA}, \ref{HB}, \ref{HC} and \ref{HD} hold. Let $N\in \NN$ satisfying
$\epsilon_{N}\leq \frac{1}{2L_{r}}(1\wedge \frac{1}{\delta})$. Then, one has 
\begin{eqnarray}
\big| \overline{\mathbfcal{H}}_{\mathcal{M}}(\mathbf{x},\mathbf{y}) -\overline{\mathbfcal{H}}_{\mathcal{M},N}(\mathbf{x},\mathbf{y}) \big| \leq K_{1} \epsilon_{N} \sum_{k=0}^{N_{0}-1} \| \mathfrak{B}^{k} \mathbf{H}\|_{\mathbb{L}(\mathcal{P}(\mathbb{X})\times\mathbb{Y})}.
\end{eqnarray}
\end{theorem}
\textbf{Proof:}
First let us show by induction that
\begin{eqnarray}
\big| \mathfrak{B}^{k}f(\theta,y) - \mathfrak{B}^{k}_{N}f(\theta,y)\big| \leq  K_1\epsilon_{N} \sum_{j=0}^{k-1} \| \mathfrak{B}^{j} f \|_{\mathbb{L}(\mathcal{P}(\mathbb{X})\times\mathbb{Y})}
\label{Diff-k-B-B-N}
\end{eqnarray}
for any $(\theta,y)\in \mathcal{P}(\mathbb{X})\times\mathbb{Y}$, $f\in \mathbb{L}(\mathcal{P}(\mathbb{X})\times\mathbb{Y})$ and $k\in \llbracket 1 ; N_{0}\rrbracket$.
From Proposition \ref{Prop-Diff-B-B-N}, the claim is true for $k=1$. Now, assume that equation (\ref{Diff-k-B-B-N}) holds for $k\in \llbracket 1 ; N_{0}-1\rrbracket$.
Then, 
\begin{align}
\big| \mathfrak{B}^{k+1}f(\theta,y) - \mathfrak{B}^{k+1}_{N}f(\theta,y)\big|
& \leq \big| \mathfrak{B}(\mathfrak{B}^{k}f)(\theta,y) - \mathfrak{B}_{N}(\mathfrak{B}^{k}f)(\theta,y)\big| \nonumber \\
& \phantom{\leq }+ \big| \mathfrak{B}_{N}(\mathfrak{B}^{k}f)(\theta,y) - \mathfrak{B}_{N}(\mathfrak{B}_{N}^{k}f)(\theta,y)\big|.
\label{Rec1-Diff-k-B-B-N}
\end{align}
From equation (\ref{Diff-B-B-N}), Lemma~\ref{B Lip} and recalling the definition of $\mathfrak{B}_{N}$ (see equation (\ref{Def-B-N})) we get
\begin{align}
\big| \mathfrak{B}^{k+1}f(\theta,y) - \mathfrak{B}^{k+1}_{N}f(\theta,y)\big|
& \leq  \|\mathfrak{B}^{k}f\|_{\mathbb{L}(\mathcal{P}(\mathbb{X})\times\mathbb{Y})}  K_{1} \epsilon_{N}  %\nonumber \\
%& \phantom{\leq }
+ \big| S_{N}(\mathfrak{B}^{k}f)(\theta,y) - S_{N}(\mathfrak{B}_{N}^{k}f)(\theta,y)\big|.
\end{align}
Now, combining (\ref{Def-S-N}) and the induction hypothesis we have
\begin{align}
\big| S_{N}(\mathfrak{B}^{k}f)(\theta,y) - S_{N}(\mathfrak{B}_{N}^{k}f)(\theta,y)\big|
& \leq \int_{\mathbb{Y}} \Big| \mathfrak{B}^{k}f(\Phi_{N}(v,\theta,y),v)-\mathfrak{B}_{N}^{k}f(\Phi_{N}(v,\theta,y),v) \Big| R_{N}(\mathbb{X},dv|\theta,y)
\nonumber \\
& \leq K_{1} \epsilon_{N} \sum_{j=0}^{k-1} \| \mathfrak{B}^{j} f \|_{\mathbb{L}(\mathcal{P}(\mathbb{X})},
\label{Rec2-Diff-k-B-B-N}
\end{align}
and so from equations (\ref{Rec1-Diff-k-B-B-N})-(\ref{Rec2-Diff-k-B-B-N}) we obtain that (\ref{Diff-k-B-B-N}) holds for any $k\in \llbracket 1 ; N_{0}\rrbracket$.
Finally, recalling that $\overline{\mathbfcal{H}}_{\mathcal{M}}(\mathbf{x},\mathbf{y}) = \mathfrak{B}^{N_{0}}\mathbf{H}(\delta_{\mathbf{x}},\mathbf{y})$ and
$\overline{\mathbfcal{H}}_{\mathcal{M},N}(\mathbf{x},\mathbf{y})= \mathfrak{B}_{N}^{N_{0}}\mathbf{H}(\delta_{\mathbf{x}},\mathbf{y})$, we obtain the result by applying (\ref{Diff-k-B-B-N}) to $f=\mathbf{H}$ with $k=N_{0}$ since $\mathbf{H}\in \mathbb{L}(\mathcal{P}(\mathbb{X})\times\mathbb{Y})$ by Assumption \ref{HD}.
\hfill $\Box$

%%%%%%%%%%%%%%%%%%%%%%%%%%%%%%%
\subsection{Second approximation}
%%%%%%%%%%%%%%%%%%%%%%%%%%%%%%%
\label{approx-2}
The approximate value function  $\overline{\mathbfcal{H}}_{\mathcal{M},N}(\mathbf{x},\mathbf{y})= \mathfrak{B}_{N}^{N_{0}}\mathbf{H}(\delta_{\mathbf{x}},\mathbf{y})$ is not directly computable as it involves a recursion of functions defined on the continuous space $\mathcal{P}(\Gamma^{N}_{X})\times\mathbb{Y}$. In order to obtain a numerically tractable recursion, one additional discretization procedure is required. We first introduce the Markov chain $\Psi^{N}$ with transition kernel $S_N$, then rewrite the iteration of the Bellman operators $\mathfrak{B}_{N}$ in terms of conditional expectations involving this chain, and finally propose and approximation of the latter conditional expectations based on the quantization of the chain $\Psi^{N}$.  Following the idea of \cite{pham05}, instead of discretizing the two coordinates (filter and observations) separately, we discretize them jointly exploiting the Markov property of $\Psi^{N}$.\\

Let us denote by $\{\Psi^{N}_{t}\}_{t\in \llbracket 0 ; N_{0}\rrbracket}$ the Markov chain with transition kernel $S_{N}$ and initial distribution $(\delta_{\mathbf{x}},\mathbf{y})$. 
By definition of $S_{N}$ we have that $S_{N}(\mathcal{P}(\Gamma^{N}_{X})\times\mathbb{Y} |\theta,y)=1$ for any $(\theta,y)\in \mathcal{P}(\mathbb{X})\times\mathbb{Y}$ and that
$\delta_{\mathbf{x}}\in \mathcal{P}(\Gamma^{N}_{X})$.
Moreover, it is clear that $\mathcal{P}(\Gamma^{N}_{X})$ can be identified with the $N$-simplex in $\RR^{N}$ denoted by
$\mathbb{S}^{N}$.
Therefore, by a slight abuse of notation %there is no loss of generality to consider that 
we will consider from now on that the state space of the Markov chain $\{\Psi^{N}_{t}\}_{t\in \llbracket 0 ; N_{0}\rrbracket}$ is given by
$\mathbb{S}^{N}\times \mathbb{Y} \subset\RR^{N+n}$. For notational convenience, the stochastic kernel associated with $\{\Psi^{N}_{t}\}_{t\in \llbracket 0 ; N_{0}\rrbracket}$ will still be denoted by $S_{N}$.
Thus, our aim is now to rewrite the Bellman operator $\mathfrak{B}_{N}$  in terms of conditional expectations involving $\{\Psi^{N}_{t}\}_{t\in \llbracket 0 ; N_{0}\rrbracket}$, discretize this Markov chain using optimal quantization and see how the approximation error is propagated through the dynamic programming recursion.\\

First, we rewrite the dynamic programming recursion on functions (\ref{value:N}) as a recursion involving conditional expectations.
By a slight abuse of notation, we write
\begin{align}
\label{Def-H-ext}
\mathbf{H}(\psi)=\sum_{j=1}^{N} \gamma_{j} \mathbf{H}(x^{j}_{N},y),
\end{align}
for $\psi=(\gamma,y) \in \mathbb{S}^{N}\times \mathbb{Y}$.
Define recursively the sequence of real-valued functions $\{V^{N}_{t}\}_{t\in \llbracket 0 ; N_{0}\rrbracket}$ on $\mathbb{S}^{N}\times \mathbb{Y}$ by
\begin{align}
\label{Def-VN-t}
V^{N}_{t}(\psi) & = \max \Big\{ \mathbf{H}(\psi),\mathbb{E} \big[V_{t+1}(\Psi^{N}_{t+1}) \big| \Psi^{N}_{t}=\psi \big] \Big\} 
\end{align}
for $t\in \llbracket 0 ; N_{0}-1\rrbracket$ and  $V^{N}_{N_{0}}(\psi) =\mathbf{H}(\psi)$ for $\psi\in \mathbb{S}^{N}\times \mathbb{Y}$.
Note that these dynamic programming equations now go backward in time, with an initialisation at the terminal time $N_0$.
By definition of the operator $\mathfrak{B}_{N}$ (see equation \ref{Def-B-N}), we have clearly
$V^{N}_{0}(\Psi^{N}_{0})=\mathfrak{B}_{N}^{N_{0}}\mathbf{H}(\delta_{\mathbf{x}},\mathbf{y})$ and so by Theorem \ref{Theo-approx-1}, $V^{N}_{0}(\Psi^{N}_{0})=\overline{\mathbfcal{H}}_{\mathcal{M},N}(\mathbf{x},\mathbf{y})$. Thus one just needs to build a numerically computable approximation of function $V^{N}_{0}$.\\

Let $\{\widehat{\Psi}^{N,M}_{t}=(\widehat{\Theta}^{N,M}_{t},\widehat{Y}^{N,M}_{t})\}_{n\in \llbracket 0 ; N_{0}\rrbracket}$ be the quantization approximation of
$\{\Psi^{N}_{t}\}_{t\in \llbracket 0 ; N_{0}\rrbracket}$ defined on a probability space $(\widebar{G},\widebar{\mathcal{G}},\widebar{\mathbb{P}})$
($\widebar{\mathbb{E}}[\cdot]$ will stand for the expectation associated to $\widebar{\mathbb{P}}$).
There are several methods to get the quantization of a Markov chain such as the marginal quantization or Markovian quantization approaches.
These techniques are roughly speaking based upon the quantization of a random variable as described in section \ref{sec:quantize}.
We do not want to go into the details of these different approaches. A rather complete exposition of this subject can be found in \cite{bally05,pages04}.
We write $\Gamma_{\Psi^{N}_{t}}^{M}$ for the grid of $M$ points used to quantize $\Psi^{N}_{t}$ and $\| \Psi^{N}_{t}-\widehat{\Psi}^{N,M}_{t}\|_{2}$
for the $L_{2}$-quantization error between $\Psi^{N}_{t}$ and $\widehat{\Psi}^{N,M}_{t}$ under $\widebar{\mathbb{P}}$.
Define recursively the sequence of real-valued functions $\{\widehat{V}^{N,M}_{t}\}_{t\in \llbracket 0 ; N_{0}\rrbracket}$ by
\begin{align*}
\widehat{V}^{N,M}_{t}(\widehat{\psi}) & = \max \Big\{ \mathbf{H}(\widehat{\psi}),
\widebar{\mathbb{E}} \big[\widehat{V}^{N,M}_{t+1}(\widehat{\Psi}^{N}_{t+1}) \big| \widehat{\Psi}^{N}_{t}=\widehat{\psi} \big] \Big\} ,
\end{align*}
for any $\widehat{\psi} \in \Gamma_{\Psi^{N}_{t}}^{M}$, $t\in \llbracket 0 ; N_{0}-1\rrbracket$ and
$\widehat{V}^{N,M}_{N_{0}}(\widehat{\psi}) = \mathbf{H}(\widehat{\psi})$ for $\widehat{\psi} \in \Gamma_{\Psi^{N}_{N_{0}}}^{M}$.
As $\{\widehat{\Psi}^{N,M}_{t}\}$ is now a (inhomogeneous) Markov chain on a finite state space, the conditional expectations above are just weighted sums and can be computed numerically. Before stating the main result of this section regarding the convergence of $\widehat{V}^{N,M}_{t}$ to $\widehat{V}^{N}_{t}$, we need additional technical results on the Lipschitz regularity of $V^N$ and $V^{N,M}$.

\begin{lemma}
\label{Lip-S-N}
Suppose Assumptions \ref{HA}, \ref{HB}, \ref{HC} and \ref{HD} hold. Let $N\in \NN$ satisfying
$\epsilon_{N}\leq \frac{1}{2L_{r}}(1\wedge \frac{1}{\delta})$. Then $V^{N}_{t}\in \mathbb{L}(\mathbb{S}^{N}\times \mathbb{Y})$ and
$\|V^{N}_{t}\|_{sup}\leq \|\mathbf{H}\|_{sup}$ for $t\in \llbracket 0 ; N_{0}\rrbracket$. Moreover, one has
\begin{align}
L_{V^{N}_{t}} \leq 4\sqrt{N} \Big[(1+2\widebar{r}) \|\mathbf{H}\|_{sup}
+2\widebar{r}\delta (1+2\widebar{r}\delta)  L_{V^{N}_{t+1}}\Big] (\widebar{r}+L_{r}) +2\sqrt{N} (\|\mathbf{H}\|_{sup}+L_{\mathbf{H}})
\label{eq-Lip-S-N}
\end{align}
for $t\in \llbracket 0 ; N_{0}-1\rrbracket$ and $L_{V^{N}_{N_{0}}}\leq 2\sqrt{N} (\|\mathbf{H}\|_{sup}+L_{\mathbf{H}})$.
\end{lemma}
\textbf{Proof:}
According to equation (\ref{Def-H-ext}), it is clear that $\|V^{N}_{N_{0}}\|_{sup}\leq \|\mathbf{H}\|_{sup}=\sup_{(x,y)\in \mathbb{X}\times \mathbb{Y}} |\mathbf{H}(x,y)|$.
Moreover, for $\psi=(\gamma,y)$ and $\psi'=(\gamma',y')$ in $\mathbb{S}^{N}\times \mathbb{Y}$
\begin{align*}
|V^{N}_{N_{0}}(\psi)-V^{N}_{N_{0}}(\psi')| & %= \big| \sum_{j=0}^{N} \gamma_{j} \mathbf{H}(x^{j}_{N},y) - \sum_{j=0}^{N} \gamma_{j}' \mathbf{H}(x^{j}_{N},y') \big|
\leq (\|\mathbf{H}\|_{sup}+L_{\mathbf{H}}) \Big[ \sum_{j=1}^{N}|\gamma_{j}-\gamma_{j}'| + \big| y-y'\big| \Big]
\nonumber \\
&  \leq \|\mathbf{H}\|_{\mathbb{L}(\mathbb{S}^{N}\times \mathbb{Y})}\Big[ \sqrt{N} |\gamma-\gamma'| + \big| y-y'\big| \Big]
\leq 2\sqrt{N} \|\mathbf{H}\|_{\mathbb{L}(\mathbb{S}^{N}\times \mathbb{Y})} |\psi-\psi'|,
\end{align*}
giving the Lipschitz constant of $V^{N}_{N_{0}}$.
\newline
Now, by a slight abuse of notation, $\Phi_{N}(Y,\gamma,y)$ is identified with the vector in $\mathbb{S}^{N}$ which the $j^{th}$ component is given by
$\Phi_{N}(Y,\gamma,y)(x_{j}^{N})$ and $r(X_{N},Y,\gamma,y)$
(respectively, $r(\lambda_{N},\nu,\gamma,y)$) denotes $\ds \sum_{j=1}^{N}\gamma_{j}r(X_{N},Y,x_{j}^{N},y)$
(respectively, $\ds \sum_{j=1}^{N}\gamma_{j}r(\lambda_{N},\nu,x_{j}^{N},y)$).
\newline
Consider $g\in \mathbb{L}(\mathbb{S}^N \times \mathbb{Y})$ and
$\psi=(\gamma,y)$, $\psi'=(\gamma',y')$ in $\mathbb{S}^N \times \mathbb{Y}$. For any $t\in \llbracket 0 ; N_{0}-1\rrbracket$, we have
\begin{align}
\Big| \mathbb{E} & \big[g(\Psi^{N}_{t+1}) \big| \Psi^{N}_{t}=\psi \big] - \mathbb{E}^{N} \big[g(\Psi^{N}_{t+1}) \big| \Psi^{N}_{t}=\psi \big] \Big|
=\big| S_{N}g(\gamma,y)-  S_{N}g(\gamma',y') \Big|
\nonumber \\
& \leq \mathbb{E} \Big[ \Big| g(\Phi_{N}(Y,\gamma,y),Y) \frac{r(\lambda_{N},Y,\gamma,y)}{r(\lambda_{N},\nu,\gamma,y)}
- g(\Phi_{N}(Y,\gamma',y'),Y) \frac{r(\lambda_{N},Y,\gamma',y')}{r(\lambda_{N},\nu,\gamma',y')} \Big| \Big]
\nonumber \\
& \leq \mathbb{E} \Big[  \frac{|g(\Phi_{N}(Y,\gamma,y),Y)|}{r(\lambda_{N},\nu,\gamma,y)} \big| r(\lambda_{N},\nu,\gamma,y) - r(\lambda_{N},\nu,\gamma',y') \big| \Big]
\nonumber \\
& \phantom{\leq} + \mathbb{E} \Big[ |g(\Phi_{N}(Y,\gamma,y),Y)| r(\lambda_{N},\nu,\gamma',y') \Big| \frac{1}{r(\lambda_{N},\nu,\gamma,y)}
- \frac{1}{r(\lambda_{N},\nu,\gamma',y')} \Big| \Big]  
\nonumber \\
& \phantom{\leq} + \mathbb{E} \Big[ \frac{r(\lambda_{N},Y,\gamma',y')}{r(\lambda_{N},\nu,\gamma',y')}  | g(\Phi_{N}(Y,\gamma,y),Y)  - g(\Phi_{N}(Y,\gamma',y'),Y) | \Big].
\label{eq1-Lip-S-N}
\end{align}
By using Lemma \ref{bound-1/r-N} and Assumption \ref{HB} we have
\begin{align}
\mathbb{E} \Big[  & \frac{|g(\Phi_{N}(Y,\gamma,y),Y)|}{r(\lambda_{N},\nu,\gamma,y)} \big| r(X_{N},\nu,\gamma,y) -   r(X_{N},\nu,\gamma',y') \big| \Big]
\nonumber \\
& \leq 2 \|g\|_{sup} (\widebar{r}+L_{r}) \Big[ \sum_{j=1}^{N}|\gamma_{j}-\gamma_{j}'| + \big| y-y'\big| \Big]
\leq 2 \|g\|_{sup} (\widebar{r}+L_{r}) \Big[ \sqrt{N} |\gamma-\gamma'| + \big| y-y'\big| \Big].
\label{eq2-Lip-S-N}
\end{align}
Similarly,
\begin{align}
\mathbb{E} \Big[ & |g(\Phi_{N}(Y,\gamma,y),Y)| r(\lambda_{N},\nu,\gamma',y')  \Big| \frac{1}{r(\lambda_{N},\nu,\gamma,y)} - \frac{1}{r(\lambda_{N},\nu,\gamma',y')} \Big| \Big]  
\nonumber \\
& \leq \mathbb{E} \Big[ |g(\Phi_{N}(Y,\gamma,y),Y)| \frac{r(\lambda_{N},\nu,\gamma',y')}{r(\lambda_{N},\nu,\gamma,y)r(\lambda_{N},\nu,\gamma',y')}
\Big| r(\lambda_{N},\nu,\gamma,y) - r(\lambda_{N},\nu,\gamma',y') \Big| \Big]  
\nonumber \\
& \leq  4 \|g\|_{sup}  \widebar{r} (\widebar{r}+L_{r})\Big[ \sqrt{N} |\gamma-\gamma'| + \big| y-y'\big| \Big],
\label{eq3-Lip-S-N}
\end{align}
and
\begin{align}
\mathbb{E} \Big[ \frac{r(\lambda_{N},Y,\gamma',y')}{r(\lambda_{N},\nu,\gamma',y')}  | g(\Phi_{N}(Y,\gamma,y),Y)  - & g(\Phi_{N}(Y,\gamma',y'),Y) | \Big]
\nonumber \\
& \leq 2 \widebar{r} L_{g}  \mathbb{E} \big[ \big| \Phi_{N}(Y,\gamma,y)-\Phi_{N}(Y,\gamma',y')\big| \big].
\label{eq4-Lip-S-N}
\end{align}
Moreover, from the definition of the discrete measure $\Phi_{N}$ (see equation (\ref{Def-Phi-N}))
\begin{align*}
\mathbb{E} \big[ \big| \Phi_{N}(Y,\gamma,y)-\Phi_{N}(Y,\gamma',y')\big| \big] & \leq 
%\sup_{v\in \mathbb{Y}}d_{\mathcal{P}}(\Phi_{N}(v,\gamma,y),\Phi_{N}(v,\gamma',y')) \big] & \leq \sup_{v\in \mathbb{Y}}
\mathbb{E} \Big[ \Big(\sum_{j=1}^{N}  \lambda_{N}(x^{j}_{N})^{2} \Big| \frac{r(x^{j}_{N},Y,\gamma,y)}{r(\lambda_{N},Y,\gamma,y)}
- \frac{r(x^{j}_{N},Y,\gamma',y')}{r(\lambda_{N},Y,\gamma',y')} \Big|^{2} \Big)^{1/2}\Big]  
\nonumber \\
& \leq \sup_{(u,v)\in \Gamma^{N}_{X}\times \mathbb{Y}}
\Big| \frac{r(u,v,\gamma,y)}{r(\lambda_{N},v,\gamma,y)} - \frac{r(u,v,\gamma',y')}{r(\lambda_{N},v,\gamma',y')} \Big| 
\nonumber \\
& \leq \sup_{(u,v)\in \Gamma^{N}_{X}\times \mathbb{Y}}
\Big[ \frac{1}{r(\lambda_{N},v,\gamma,y)} \Big| r(u,v,\gamma,y) - r(u,v,\gamma',y') \Big| \Big]  
\nonumber \\
& \phantom{\leq} +\sup_{(u,v)\in \Gamma^{N}_{X}\times \mathbb{Y}}
 \Big[ r(u,v,\gamma',y')  \Big| \frac{1}{r(\lambda_{N},v,\gamma,y)} - \frac{1}{r(\lambda_{N},v,\gamma',y')} \Big| \Big] ,
\end{align*}
and so, from Lemma \ref{bound-1/r-N-bis}
\begin{align}
\mathbb{E} \big[ \big| \Phi_{N}(Y,\gamma,y)-\Phi_{N}(Y,\gamma',y')\big| \big] 
& \leq 2\delta (1+2\widebar{r}\delta) \sup_{(u,v)\in \Gamma^{N}_{X}\times \mathbb{Y}}
 \Big| r(u,v,\theta,y) - r(u,v,\theta',y') \Big| 
 \nonumber \\
&\leq  2\delta (1+2\widebar{r}\delta) (\widebar{r}+L_{r}) \Big[ \sqrt{N} |\gamma-\gamma'| + \big| y-y'\big| \Big].
\label{eq5-Lip-S-N}
\end{align}
Combining equations (\ref{eq1-Lip-S-N})-(\ref{eq5-Lip-S-N}), we obtain
\begin{align*}
\Big| \mathbb{E} & \big[g(\Psi^{N}_{t+1}) \big| \Psi^{N}_{t}=\psi \big] - \mathbb{E}^{N} \big[g(\Psi^{N}_{t+1}) \big| \Psi^{N}_{t}=\psi \big] \Big|
=\big| S_{N}g(\gamma,y)-  S_{N}g(\gamma',y') \Big|
\nonumber \\
& \leq 4 \sqrt{N} \Big[(1+2\widebar{r}) \|g\|_{sup} +2\widebar{r}\delta (1+2\widebar{r}\delta) L_{g} \Big] (\widebar{r}+L_{r})  \big| \psi-\psi'\big|.
\end{align*}
Finally, by using the definition of $V^{N}_{t}$ (see equation (\ref{Def-VN-t})), we get (\ref{eq-Lip-S-N}) showing the result.
\hfill $\Box$\\

We now state and prove the main result of this section.
\begin{theorem}
\label{Theo-approx-2}
Suppose Assumptions \ref{HA}, \ref{HB}, \ref{HC} and \ref{HD} hold. Let $N\in \NN$ satisfying
$\epsilon_{N}\leq \frac{1}{2L_{r}}(1\wedge \frac{1}{\delta})$. Then
\begin{eqnarray*}
\big| \overline{\mathbfcal{H}}_{\mathcal{M},N}(\mathbf{x},\mathbf{y}) -\widehat{V}^{N,M}_{0}(\widehat{\Psi}^{N,M}_{0}) \big| \leq \sum_{t=0}^{N_{0}} L_{V^{N}_{t}}  \| \Psi^{N}_{t}-\widehat{\Psi}^{N,M}_{t}\|_{2}.
\end{eqnarray*}
\end{theorem}
\textbf{Proof:} The proof of this result is based on Theorem 2 in \cite{bally05}. The main difference is that in our setting, the transition kernel of Markov chain $\{\Psi^{N}_{t}\}_{t\in \llbracket 0 ; N_{0}\rrbracket}$ is not \textit{K-Lipschitz} in the sense of the definition (2.13) in \cite{bally05}. However, the main arguments of the proof of Theorem 2 in \cite{bally05} can still be applied to show that
\begin{align*}
\big\|V^{N}_{t}(\Psi^{N}_{t})-\widehat{V}^{N,M}_{t}(\widehat{\Psi}^{N,M}_{t})\big\|_{2} \leq L_{V^{N}_{t}}  \big\| \Psi^{N}_{t}-\widehat{\Psi}^{N,M}_{t}\big\|_{2} + \big\|V^{N}_{t+1}(\Psi^{N}_{t+1})-\widehat{V}^{N,M}_{t+1}(\widehat{\Psi}^{N,M}_{t+1})\big\|_{2},
\end{align*}
for $t\in \llbracket 0 ; N_{0}-1\rrbracket$
and 
\begin{align*}
\big\|V^{N}_{N_{0}}(\Psi^{N}_{_{0}})-\widehat{V}^{N,M}_{N_{0}}(\widehat{\Psi}^{N,M}_{N_{0}})\big\|_{2} \leq L_{V^{N}_{N_{0}}}  \big\| \Psi^{N}_{N_{0}}-\widehat{\Psi}^{N,M}_{N_{0}}\big\|_{2} ,
\end{align*}
where $L_{V^{N}_{t}}$ are given in Lemma \ref{Lip-S-N}.
This implies that
\begin{eqnarray*}
\big| V^{N}_{0}(\Psi^{N}_{0}) -\widehat{V}^{N,M}_{0}(\widehat{\Psi}^{N,M}_{0}) \big| \leq \sum_{t=0}^{N_{0}} L_{V^{N}_{t}}  \| \Psi^{N}_{t}-\widehat{\Psi}^{N,M}_{t}\|_{2}.
\end{eqnarray*}
Moreover, one has $V^{N}_{0}(\Psi^{N}_{0})=\overline{\mathbfcal{H}}_{\mathcal{M},N}(\mathbf{x},\mathbf{y})$
giving the result.
\hfill $\Box$\\

Gathering together our three main results Theorems~\ref{Equivalence-Costs}, \ref{Theo-approx-1}, and \ref{Theo-approx-2}, we obtain that the fully computable expression $\widehat{V}^{N,M}_{0}(\widehat{\Psi}^{N,M}_{0})$ is an approximation of our initial value function of interest $\overline{\mathbfcal{H}}(\mathbf{x},\mathbf{y})$ with an error bound of
\begin{equation*}
\big|\overline{\mathbfcal{H}}(\mathbf{x},\mathbf{y})-\widehat{V}^{N,M}_{0}(\widehat{\Psi}^{N,M}_{0})\big|
\leq K_{1} \epsilon_{N} \sum_{k=0}^{N_{0}-1} \| \mathfrak{B}^{k} \mathbf{H}\|_{\mathbb{L}(\mathcal{P}(\mathbb{X})\times\mathbb{Y})}
+\sum_{t=0}^{N_{0}} L_{V^{N}_{t}}  \| \Psi^{N}_{t}-\widehat{\Psi}^{N,M}_{t}\|_{2},
\end{equation*}
that goes to zero as the number of points in the quantization grids goes to infinity.

\section{Numerical example}\label{sec:example}

In this section, we present a numerical example to illustrate our approximation results. It is adapted from the control of water tank problems which can be found in \cite[section 1.3]{hernandez89}. Such applications are essential in regions under high water stress.

Consider a water tank which capacity $K > 0$ is finite. It is filled with a random amount of rainfall each time it rains. However, the water level is only known through noisy measurements. One wants to cover the tank when the volume of water is closest to some value $\alpha \in (0;K)$.
Let us model this situation with a $[0;K]^2$-valued finite-horizon Markov chain $(\mathbfcal{\tilde{X}}_t,\mathbfcal{\tilde{Y}}_t)_{t \in \llbracket 0;N_0 \rrbracket}$, where $(\mathbfcal{\tilde{X}}_t)$ represents the sequence of water volumes contained in the tank and $(\mathbfcal{\tilde{Y}}_t)$ symbolizes the measurements of $(\mathbfcal{\tilde{X}}_t)$. We suppose that the dynamics of the Markov chain is given by
\begin{equation*}
\left\{
\begin{array}{l}
\mathbfcal{\tilde{X}}_{t+1} = \min\big\{(\mathbfcal{\tilde{X}}_t+ \xi_t)_+; K\big\}\\
\mathbfcal{\tilde{Y}}_{t+1} = \min\big\{(\mathbfcal{\tilde{X}}_{t+1} + \psi_t)_+ ; K\big\} 
\end{array}
\right.
\end{equation*}
where $x_+$ stands for the positive part of a real number $x$, and $(\xi_t)$ and $(\psi_t)$ are i.i.d. random variables with respective densities $f$ on $\mathbb{R}_+$ and $g$ on $\mathbb{R}$. Let us denote respectively $F$ and $G$ the cumulative distribution functions associated to $f$ and $g$.
Let $B,C \in \mathcal{B}([0;K])$. The cost function is $\mathbf{\tilde{H}}(x,y) = K - |x - \alpha|$ so that the process is optimally stopped when the (unobserved) component $\mathbfcal{\tilde{X}}_t$ is close or equal to $\alpha$. The transition law of this process is
\begin{equation*}
\tilde{R}(B \times C |x,y) = \delta_0(B)F(-x)\mathfrak{M}_1(C) + \int_B f(\xi - x)\mathfrak{M}_2(C,\xi)d\xi + \delta_K(B)\mathfrak{M}_3(C)
\end{equation*}
where
\begin{eqnarray*}
\mathfrak{M}_1(C) &=& \delta_0(C)G(0) + \int_C g(\psi)d\psi + \delta_K(C)(1 - G(K)),\\
\mathfrak{M}_2(C,\xi) &=& \delta_0(C)G(-\xi) + \int_C g(\psi - \xi)d\psi + \delta_K(C)(1 - G(K - \xi)),\\
\mathfrak{M}_3(C) &=& \delta_0(C)G(-K) + \int_C g(\psi - K)d\psi + \delta_K(C)(1 - G(0)).
\end{eqnarray*}
Assumption \ref{HB} does not hold when $[0;K]$ is endowed with the usual Euclidian norm because the points 0 and $K$ have a nonzero weight. 
Thus we change the topology to isolate these two points by adding an additional dimension to the process.\\

Consider the process $(\mathbfcal{X}^1_t,\mathbfcal{X}_t^2,\mathbfcal{Y}_t)_{t \in \llbracket 0;N_0 \rrbracket}$, where $\mathbfcal{X}_t^1 = \mathbfcal{\tilde{X}}_t$, $\mathbfcal{Y}_t = \mathbfcal{\tilde{Y}}_t$ and the dynamics of $\mathbfcal{X}_t^2$ is
\begin{equation*}
\mathbfcal{X}^2_{t+1} = I_{\{\mathbfcal{X}^1_{t+1} = K\}} - I_{\{\mathbfcal{X}^1_{t+1} = 0\}}.
\end{equation*}
So, the unobservable state space is $\mathbb{X} = ((0;K) \times \{0\}) \cup \{(0,-1)\} \cup \{(K,1)\}$. The observable state space is $\mathbb{Y} = [0;K]$. %
Let  $\mathbf{H}(x_1,x_2,y) = K - |x_1 - \alpha|$ be the performance function. 
One may now write the transition law $R$ of the process $(\mathbfcal{X}^1_t,\mathbfcal{X}_t^2,\mathbfcal{Y}_t)_{t \in \llbracket 0;N_0 \rrbracket}$ as
$$R(du_1,du_2,dv | x_1,x_2,y) = r(u_1,u_2,v,x_1,x_2,y)\lambda(du_1,du_2)\nu(dv),$$
where $r$ : $(\mathbb{X} \times \mathbb{Y})^2 \to \mathbb{R}_+$ is defined by
\begin{eqnarray*}
\lefteqn{r(u_1,u_2,v,x_1,x_2,y) }\\
&=&4I_{\{(0,-1)\}}(u_1,u_2)F(-x_1)\mathfrak{m}_1(v)+ 2KI_{(0;K) \times \{0\}}(u_1,u_2)f(u_1 - x_1)\mathfrak{m}_2(u_1,v)\nonumber\\
&&+ 4I_{\{(K,1)\}}(u_1,u_2)(1 - F(K - x_1))\mathfrak{m}_3(v),
\end{eqnarray*}
with
\begin{eqnarray*}
\mathfrak{m}_1(v) &=& 4G(0)I_{\{0\}}(v) + 2KI_{(0;K)}(v)g(v) + 4(1-G(K))I_{\{K\}}(v),\\
\mathfrak{m}_2(u,v) &=& 4G(-u)I_{\{0\}}(v) + 2KI_{(0;K)}(v)g(v - u) + 4(1 - G(K - u))I_{\{K\}}(v),\\
\mathfrak{m}_3(v) &=& 4G(-K)I_{\{0\}}(v) + 2KI_{(0;K)}(v)g(v - K) + 4(1 - G(0))I_{\{K\}}(v),
\end{eqnarray*}
and
\begin{eqnarray*}
\lambda(du_1,du_2) &=& \dfrac{\delta_0(du_1)\delta_{-1}(du_2)}{4} + \dfrac{\mu(du_1)\delta_0(du_2)}{2K} + \dfrac{\delta_K(du_1)\delta_1(du_2)}{4},\\
\nu(dv) &=& \dfrac{\delta_0(dv)}{4} + \dfrac{\mu(dv)}{2K} + \dfrac{\delta_K(dv)}{4},
\end{eqnarray*}
where $\mu$ denotes the Lebesgue measure. One may note that neither $\lambda$ nor $\nu$ are absolutely continuous with respect to the Lebesgue measure on $[0;K]$. Thus, this model does not satisfy the assumptions of \cite{YZ13,Zhou13,ZFM10}. However, our Assumption \ref{HA} is clearly satisfied.\\

Assume that $f$ is Lipschitz-continuous on $[0;K]$ with constant $L_f$ and $g$ is positive and Lipschitz-continuous on $[-K;K]$ with constant $L_g$ (e.g. if $f$ is an exponential density function and $g$ a centered Gaussian density function, these hold). Therefore, they are both bounded above on these intervals, respectively by $\|f\|_{sup}$ and $\|g\|_{sup}$. Straightforward calculations show that assumptions B and D hold with the following constants
\begin{eqnarray*}
\overline{r} &=& (8 + 2K\|f\|_{sup})(8 + 2K\|g\|_{sup})\\
L_r &=& \max\Big\{\mathfrak{a}, \mathfrak{b}, \dfrac{8(8 + 2K\|g\|_{sup})}{K+2}, (8 + 2K\|g\|_{sup})(8\|f\|_{sup} + 2KL_f) \Big\},\\
\|\mathbf{H}\|_{\mathbb{L}(\mathbb{X}\times \mathbb{Y})}&\leq&K+1,
\end{eqnarray*}
where
$$\mathfrak{a} = 2K(L_f(8 + 2K\|g\|_{sup}) + \|f\|_{sup}(8\|g\|_{sup} + 2KL_g))$$
and
$$\mathfrak{b} = \max\big\{2K\|f\|_{sup}(8\|g\|_{sup} + 2KL_g) ; (4 + 2K\|f\|_{sup} + 2KL_f)(8 + 2K\|g\|_{sup}) \big\}.$$
Assumption C requires that the density $g$ be bounded from below by some positive number $\mathfrak{g}$ on $[-K;K]$. Thus, one has
\begin{equation*}
\mathfrak{m}_3(v) \geq \min\{4G(-K) ; 2K\mathfrak{g} ; 4(1 - G(0)) \}.
\end{equation*}
As
$$r(\lambda,v,x_1,x_2,y) = \int_{\mathbb{X}} r(u_1,u_2,v,x_1,x_2,y)\lambda(du_1,du_2) \geq (1 - F(K - x_1))\mathfrak{m}_3(v),$$
let us suppose that $F(K) < 1$, $G(0) < 1$ and $G(-K) > 0$. These are verified by exponential and centered Gaussian density functions as above, for instance. One then deduces that
$$r(\lambda,v,x_1,x_2,y) \geq (1 - F(K))\min\big(4G(-K) ; 2K\mathfrak{g} ; 4(1 - G(0)) \big) > 0.$$
for all $v,y \in \mathbb{Y}$ and $(x_1,x_2) \in \mathbb{X}$. Therefore, this shows that assumption C holds.\\

For our numerical experimentations, we chose $N_0 = 10$, $K = 1$, $\alpha = 0.5$, and the initial state $(\mathbfcal{X}^1_0,\mathbfcal{X}^2_0,\mathbfcal{Y}_0) = (0,-1,0)$. We suppose that the $\xi_t$ are exponentially distributed with parameter 5. We suppose that the $\psi_t$ are normally distributed with mean 0 and standard deviation $0.03$.
Following the method developed in this paper, we have performed the two quantizations by using the competitive learning vector quantization algorithm (see section 2.2 of \cite{PPP04}). Table \ref{numerical} displays the approximation $\widehat{V}^{N,M}_{0}(\widehat{\Psi}^{N,M}_{0})$ of the value function at $(0,-1,0)$ according to the numbers $N$ and $M$ of points in the quantization grids. The exact value function is not known, but as expected one sees that our approximation is close to the optimal performance of $1$.

\begin{table}[t]
\begin{center}
\begin{tabular}{c | c ccccc}
 $M$ & $N = 12$ & $N = 25$ &$N=50$&$N=100$\\
\hline
%62 & 0.9299 & 0.946&&\\
125 & 0.9323 & 0.9534&&\\
250 & 0.9381 & 0.9579&&\\
500 & 0.9392 & 0.9577&&\\
1000 & 0.9404 & 0.9574&&\\
10000 & 0.9416 & 0.9578&0.9686&0.9771
\end{tabular}
\end{center}
\caption{Approximation $\widehat{V}^{N,M}_{0}(\widehat{\Psi}^{N,M}_{0})$ of the optimal value according to $M$ and $N$}
\label{numerical}
\end{table}

%%%%%%%%%%%%%%%%%%%%%%%%%%%%%%%%%%%%%%%%%%%%%%%%%%%%%%%%%%%%%
\appendix
\makeatletter
\def\@seccntformat#1{Appendix~\csname the#1\endcsname}
\makeatother
%%%%%%%%%%%%%%%%%%%%%%%%%%%%%%%%%%%%%%%%%%%%%%%%%%%%%%%%%%%%%
\section{: Proof of Theorem \ref{Equivalence-Costs}}
\label{aux-part-obs}
In order to prove Theorem \ref{Equivalence-Costs}, we need to introduce a new auxiliary control model $\mathbfcal{M}$ given by
the five-tuple $\big(\mathbb{F}, \mathbb{A},T,H,h\big)$ where
\begin{enumerate}
\item [(a)] the state space is $\mathbb{F}=\mathbb{X}\times\mathbb{Y}\times\{0,1\}$,
\item [(b)] the action space is $\mathbb{A}=\{0,1\}$,
\item [(c)] the transition probability function is given the stochastic kernel $T$ on $\mathbb{F}$ given $\mathbb{F}\times\mathbb{A}$ defined by
$T(B\times C |x,y,z,a)=R(B\times C |x,y) \big[ \delta_{z}(D) \I{a=0}+ \delta_{1}(D) \I{a=1} \big]$ for any
any $B\in \mathcal{B}(\mathbb{X})$, $C\in \mathcal{B}(\mathbb{Y})$, $D\subset\{0,1\}$ and $(x,y,z,a)\in \mathbb{F}\times\mathbb{A}$,
\item [(d)]  the cost-per-stage $H$ and the terminal cost $h$.
\end{enumerate}
Define $\mathbf{\Omega}=\mathbb{F}^{N_{0}+1}$ and $\mathbfcal{F}$ its associated product $\sigma$-algebra.
Introduce the coordinate projections $\mathbf{X}_{t}$ (respectively $\mathbf{Y}_{t}$, and $\mathbf{Z}_{t}$) from
$\mathbf{\Omega}$ to the set $\mathbb{X}$ (respectively $\mathbb{Y}$, and $\{0,1\}$).
Consider an arbitrary policy $\pi\in\Pi^{o}$. Define recursively the action process $\{\mathbf{A}_{t}\}_{t\in \llbracket 0 ; N_{0}-1 \rrbracket}$ by
$\mathbf{A}_{t}=\pi_{t}(\mathbf{Y}_{0},\mathbf{Z}_{0},\mathbf{A}_{0},\ldots, \mathbf{Y}_{t-1},\mathbf{Z}_{t-1},\mathbf{A}_{t-1},\mathbf{Y}_{t},\mathbf{Z}_{t})$
for $t\in \llbracket 1 ; N_{0}-1 \rrbracket$ and $\mathbf{A}_{0}=\pi_{0}(\mathbf{Y}_{0},\mathbf{Z}_{0})$.
Define the filtration $\{\mathbfcal{F}_{t}\}_{t\in \llbracket 0 ; N_{0} \rrbracket}$ by
%$\mathbfcal{O}_{0}=\sigma\{\mathbf{Y}_{0},\mathbf{Z}_{0}\}$,
$\mathbfcal{F}_{t}=\sigma\{\mathbf{X}_{0},\mathbf{Y}_{0},\mathbf{Z}_{0},\ldots,\mathbf{X}_{t},\mathbf{Y}_{t},\mathbf{Z}_{t}\}$ for $t\in  \llbracket 0 ; N_{0} \rrbracket$.
%and
%$$\mathbfcal{O}_{t}=\sigma\{\mathbf{Y}_{0},\mathbf{Z}_{0},\mathbf{A}_{0},\ldots,\mathbf{Y}_{t-1},\mathbf{Z}_{t-1},\mathbf{A}_{t-1},\mathbf{Y}_{t},\mathbf{Z}_{t}\}$$ for $t\in  \llbracket 1 ; N_{0} \rrbracket$
According to \cite{bauerle11,hernandez96}, there exists a probability measure $\mathbf{P}^{\pi}_{(\mathbf{x},\mathbf{y})}$
on $\big( \mathbf{\Omega},\mathbfcal{F} \big)$ satisfying
\begin{itemize}
\item[i)] $\mathbf{P}^{\pi}_{(\mathbf{x},\mathbf{y})}\big((\mathbf{X}_{0},\mathbf{Y}_{0},\mathbf{Z}_{0})\in B\times C \times D\big)=
\delta_{(\mathbf{x},\mathbf{y})}(B\times C) \delta_{0}(D)$,
\item[ii)] $\mathbf{P}^{\pi}_{(\mathbf{x},\mathbf{y})}\big((\mathbf{X}_{t+1},\mathbf{Y}_{t+1},\mathbf{Z}_{t+1})\in B\times C\times D |\mathbfcal{F}_{t} \big) =T(B\times C\times D|\mathbf{X}_{t},\mathbf{Y}_{t},\mathbf{Z}_{t},\mathbf{A}_{t})$,
\end{itemize}
for $t\in \llbracket 0 ; N_{0}-1 \rrbracket$, $B\in \mathcal{B}(\mathbb{X})$, $C\in \mathcal{B}(\mathbb{Y})$, $D\subset\{0,1\}$.

\bigskip

\noindent
The expectation under the probability $\mathbf{P}^{\pi}_{(\mathbf{x},\mathbf{y})}$ is denoted by $\mathbf{E}^{\pi}_{(\mathbf{x},\mathbf{y})}$.
For a policy $\pi\in \Pi^{o}$, the performance criterion is given by
\begin{eqnarray}
\mathbfcal{H}_{\mathbfcal{M}}(\mathbf{x},\mathbf{y},\pi) & =  &
\mathbf{E}^{\pi}_{(\mathbf{x},\mathbf{y})} \Big[ \sum_{t=0}^{N_{0}-1}  H(\mathbf{X}_{t},\mathbf{Y}_{t},\mathbf{Z}_{t},\mathbf{A}_{t}) \Big]
+\mathbf{E}^{\pi}_{(\mathbf{x},\mathbf{y})}\big[ h(\mathbf{X}_{N_{0}},\mathbf{Y}_{N_{0}},\mathbf{Z}_{N_{0}}) \big].
\end{eqnarray}
The optimization problem we are interested in is to maximize the reward function $\mathbfcal{H}_{\mathbfcal{M}}(\mathbf{x},\mathbf{y},\pi)$ over $\Pi^{o}$ and
$\overline{\mathbfcal{H}}_{\mathbfcal{M}}(\mathbf{x},\mathbf{y})=\sup_{\pi\in \Pi^{o}} \mathbfcal{H}_{\mathbfcal{M}}(\mathbf{x},\mathbf{y},\pi)$.
We first need to prove the following technical lemma.
\begin{lemma}
\label{algebra=}
For any $t\in \llbracket 0 ; N_{0}\rrbracket$,
$$\sigma\{\mathbf{Y}_{0},\mathbf{Z}_{0},\ldots,\mathbf{Y}_{t},\mathbf{Z}_{t}\} = \sigma\{\mathbf{Y}_{0},\ldots,\mathbf{Y}_{t}\}.$$
\end{lemma}
\textbf{Proof:}
Clearly, $\mathbf{A}_{t}$ is measurable with respect to $\sigma\{\mathbf{Y}_{0},\mathbf{Z}_{0},\ldots,\mathbf{Y}_{t},\mathbf{Z}_{t}\}$
for $t\in \llbracket 0 ; N_{0}-1 \rrbracket$. Moreover, from the definition of the transition kernel $T$, we obtain that $\mathbf{Z}_{t}=I_{\{\mathbf{A}_{t-1}=1\}}+\mathbf{Z}_{t-1}I_{\{\mathbf{A}_{t-1}=0\}}$ for any $t\in \llbracket 1 ; N_{0} \rrbracket$. Recalling that $\mathbf{Z}_{0}=0$, it follows easily
$\sigma\{\mathbf{Y}_{0},\mathbf{Z}_{0},\ldots,\mathbf{Y}_{t},\mathbf{Z}_{t}\} \subset \sigma\{\mathbf{Y}_{0},\ldots,\mathbf{Y}_{t}\}$
for $t\in \llbracket 0 ; N_{0}\rrbracket$ showing the result.
\hfill $\Box$

The next result shows that the optimization problem defined through $\mathbfcal{M}$ is equivalent to the initial optimal stopping problem defined in Definition \ref{Def-opt-stop}.
\begin{proposition}
\label{Equiv-partially-complete}
The following assertions hold.
\nl
$i)$ For any control $\ell \in L$, there exist a policy $\pi\in \Pi^{o}$ such that
$$\mathbfcal{H}_{\mathbfcal{M}}(\mathbf{x},\mathbf{y},\pi)=\mathbfcal{H}(\mathbf{x},\mathbf{y},\ell).$$
$ii)$ For any policy $\pi\in \Pi^{o}$, there exist a control $\ell \in L$ such that
$$\mathbfcal{H}(\mathbf{x},\mathbf{y},\ell)=\mathbfcal{H}_{\mathbfcal{M}}(\mathbf{x},\mathbf{y},\pi).$$
\end{proposition}
\textbf{Proof:}
Regarding item $i)$, consider a control $\ell =\big( \mathbf{\Xi},\mathbfcal{G},\mathbf{Q},\{\mathbfcal{G}_{t}\}_{t\in \llbracket 0 ; N_{0} \rrbracket},\{\mathbfcal{X}_{t},\mathbfcal{Y}_{t}\}_{t\in \llbracket 0 ; N_{0} \rrbracket},\tau\big)$ in $L$.
On the probability space $\big( \mathbf{\Xi},\mathbfcal{G},\mathbf{Q}\big)$, let us define the processes
$\{\mathbfcal{A}_{t}\}_{t\in \llbracket 0 ; N_{0}-1 \rrbracket}$ and $\{\mathbfcal{Z}_{t}\}_{t\in \llbracket 0 ; N_{0} \rrbracket}$ by
$\mathbfcal{A}_{t}=I_{\{\tau\leq t\}}$ and $\mathbfcal{Z}_{t}=\mathbfcal{A}_{t-1}$ for $t\in \llbracket 1 ; N_{0} \rrbracket$ and $\mathbfcal{Z}_{0}=0$.
Introduce the filtrations $\{\mathbfcal{T}_{t}\}_{t\in \llbracket 0 ; N_{0} \rrbracket}$ by
$\mathbfcal{T}_{t}=\sigma\{\mathbfcal{X}_{0},\mathbfcal{Y}_{0},\mathbfcal{Z}_{0},\mathbfcal{A}_{0},\ldots,\mathbfcal{X}_{t},\mathbfcal{Y}_{t},\mathbfcal{Z}_{t},\mathbfcal{A}_{t}\}$
and $\{\mathbfcal{G}^{\mathbfcal{Y}}_{t} \}_{t\in \llbracket 0 ; N_{0} \rrbracket}$ by $\mathbfcal{G}^{\mathbfcal{Y}}_{t}=\sigma\{\mathbfcal{Y}_{0},\ldots,\mathbfcal{Y}_{t}\}$.
Since $\tau$ is an $\{\mathbfcal{G}^{\mathbfcal{Y}}_{t}\}_{t\in \llbracket 0 ; N_{0} \rrbracket}$-stopping time, we have $\mathbfcal{T}_{t}\subset \mathbfcal{G}_{t}$.
Moreover, $\mathbfcal{Z}_{t+1}$ is $\mathbfcal{T}_{t}$-measurable. Consequently, it is easy to show that
$$\mathbf{Q}\big((\mathbfcal{X}_{t+1},\mathbfcal{Y}_{t+1},\mathbfcal{Z}_{t+1})\in B\times C\times D | \mathbfcal{T}_{t} \big)
=I_{\{\mathbfcal{Z}_{t+1}\in D\}} R(B\times C | \mathbfcal{X}_{t},\mathbfcal{Y}_{t}).$$
We have $\{\mathbfcal{A}_{t}=1\}=\{\mathbfcal{Z}_{t+1}=1\}$ and
$\{\mathbfcal{A}_{t}=0\}=\{\mathbfcal{Z}_{t+1}=0\}\subset\{\mathbfcal{A}_{t-1}=0\}=\{\mathbfcal{Z}_{t}=0\}$, and so
\begin{align}
\mathbf{Q}\big((\mathbfcal{X}_{t+1},\mathbfcal{Y}_{t+1},\mathbfcal{Z}_{t+1})\in B\times C\times D | \mathbfcal{T}_{t} \big)
= & \big[I_{\{\mathbfcal{A}_{t}=0\}} \delta_{\{\mathbfcal{Z}_{t}\in D\}} +I_{\{\mathbfcal{A}_{t}=1\}} \delta_{1}(D)\big]
R(B\times C | \mathbfcal{X}_{t},\mathbfcal{Y}_{t}) \nonumber\\
= & T(B\times C \times D | \mathbfcal{X}_{t},\mathbfcal{Y}_{t},\mathbfcal{Z}_{t},\mathbfcal{A}_{t}).
\label{egal-eq1}
\end{align}
Now, there exists an $\mathbb{A}$-valued measurable mapping $\pi_{t}$
defined on $\mathbb{Y}^{t+1}$ satisfying $\mathbfcal{A}_{t}=\pi_{t}(\mathbfcal{Y}_{0},\ldots,\mathbfcal{Y}_{t})$ and so,
\begin{align}
\mathbf{Q}(\mathbfcal{A}_{t}\in F| \sigma\{\mathbfcal{Y}_{0},\mathbfcal{Z}_{0},\mathbfcal{A}_{0},\ldots,\mathbfcal{Y}_{t},\mathbfcal{Z}_{t}\} )=\delta_{\pi_{t}(\mathbfcal{Y}_{0},\ldots,\mathbfcal{Y}_{t})}(F),
\label{egal-eq2}
\end{align}
for any $t\in \llbracket 0 ; N_{0}-1 \rrbracket$ and $F\subset \mathbb{A}$.
Recall that 
\begin{align}
\mathbf{Q}\big((\mathbfcal{X}_{0},\mathbfcal{Y}_{0},\mathbfcal{Z}_{0})\in B\times C \times D\big)=
\delta_{(\mathbf{x},\mathbf{y})}(B\times C) \delta_{0}(D)
\label{egal-eq3}
\end{align}
for any $B\in \mathcal{B}(\mathbb{X})$, $C\in \mathcal{B}(\mathbb{Y})$, $D\subset\{0,1\}$.
Combining equations (\ref{egal-eq1})-(\ref{egal-eq3}) and by the uniqueness property in the Theorem of Ionescu-Tulcea (see, e.g. \cite[Proposition C.10]{hernandez96}), it follows that for the control policy $\pi=\{\pi_{t}\}_{t\in \llbracket 0 ; N_{0} \rrbracket}$
\begin{align}
\mathbf{Q}\big((\mathbfcal{X}_{0}, &\mathbfcal{Y}_{0},\mathbfcal{Z}_{0},\mathbfcal{A}_{0},\ldots, \mathbfcal{X}_{N_{0}-1},\mathbfcal{Y}_{N_{0}-1},\mathbfcal{Z}_{N_{0}-1},\mathbfcal{A}_{N_{0}-1},
\mathbfcal{X}_{N_{0}},\mathbfcal{Y}_{N_{0}},\mathbfcal{Z}_{N_{0}})\in H\big)
\nonumber \\
& =\mathbf{P}^{\pi}_{(\mathbf{x},\mathbf{y})}\big((\mathbf{X}_{0},\mathbf{Y}_{0},\mathbf{Z}_{0},\mathbf{A}_{0},\ldots,\mathbf{X}_{N_{0}-1},\mathbf{Y}_{N_{0}-1},\mathbf{Z}_{N_{0}-1},\mathbf{A}_{N_{0}-1},
\mathbf{X}_{N_{0}},\mathbf{Y}_{N_{0}},\mathbf{Z}_{N_{0}})\in H \big)
\label{distribution=}
\end{align}
for any $H\in \mathbfcal{F}$.

Observe that for $k\in \llbracket 0 ; N_{0}-1 \rrbracket$ we have $\{\tau=k\}=\{\mathbf{Z}_{k}=0\}\cup\{\mathbf{A}_{k}=1\}$ and $\{\tau=N_{0}\}=\{\mathbf{Z}_{N_{0}}=0\}$.
Consequently,
\begin{align*}
\mathbf{E}^{\mathbf{Q}}_{(\mathbf{x},\mathbf{y})}\big[ \mathbf{H}(\mathbfcal{X}_{\tau},\mathbfcal{Y}_{\tau})\big]
& = \sum_{t=0}^{N_{0}-1} \mathbf{E}^{\mathbf{Q}}_{(\mathbf{x},\mathbf{y})}\big[ \mathbf{H}(\mathbfcal{X}_{t},\mathbfcal{Y}_{t}) I_{\{\tau=t\}}\big]
+\mathbf{E}^{\mathbf{Q}}_{(\mathbf{x},\mathbf{y})}\big[ \mathbf{H}(\mathbfcal{X}_{N_{0}},\mathbfcal{Y}_{N_{0}}) I_{\{\tau=N_{0}\}}\big] \nonumber \\
& = \sum_{t=0}^{N_{0}-1} \mathbf{E}^{\mathbf{Q}}_{(\mathbf{x},\mathbf{y})}\big[ \mathbf{H}(\mathbfcal{X}_{t},\mathbfcal{Y}_{t}) I_{\{(\mathbfcal{Z}_{t},\mathbfcal{A}_{t})=(0,1)\}}\big]
+\mathbf{E}^{\mathbf{Q}}_{(\mathbf{x},\mathbf{y})}\big[ \mathbf{H}(\mathbfcal{X}_{N_{0}},\mathbfcal{Y}_{N_{0}}) I_{\{\mathbfcal{Z}_{N_{0}}=1\}}\big]
\end{align*}
Now, by using the definitions of $H$ and $h$ we get
\begin{align*}
\mathbf{E}^{\mathbf{Q}}_{(\mathbf{x},\mathbf{y})}\big[ \mathbf{H}(\mathbfcal{X}_{\tau},\mathbfcal{Y}_{\tau})\big]
& = \sum_{t=0}^{N_{0}-1} \mathbf{E}^{\mathbf{Q}}_{(\mathbf{x},\mathbf{y})}\big[ H(\mathbfcal{X}_{t},\mathbfcal{Y}_{t},\mathbfcal{Z}_{t},\mathbfcal{A}_{t}) \big]
+\mathbf{E}^{\mathbf{Q}}_{(\mathbf{x},\mathbf{y})}\big[ h(\mathbfcal{X}_{N_{0}},\mathbfcal{Y}_{N_{0}},\mathbfcal{Z}_{N_{0}}) \big].
\end{align*}
By using equation (\ref{distribution=}), it follows that
\begin{align*}
\mathbf{E}^{\mathbf{Q}}_{(\mathbf{x},\mathbf{y})}\big[ \mathbf{H}(\mathbfcal{X}_{\tau},\mathbfcal{Y}_{\tau})\big]
& = \sum_{t=0}^{N_{0}-1}  \mathbf{E}^{\pi}_{(\mathbf{x},\mathbf{y})} \big[ H(\mathbf{X}_{t},\mathbf{Y}_{t},\mathbf{Z}_{t},\mathbf{A}_{t}) \big]
+\mathbf{E}^{\pi}_{(\mathbf{x},\mathbf{y})}\big[ h(\mathbf{X}_{N_{0}},\mathbf{Y}_{N_{0}},\mathbf{Z}_{N_{0}}) \big],
\end{align*}
showing the first claim.

\bigskip

\noindent
Regarding  item $ii)$, let $\pi$ be a policy in $\Pi^{o}$. Then, on the probability space
$\big(\mathbf{\Omega},\mathbfcal{F},\mathbf{P}^{\pi}_{(\mathbf{x},\mathbf{y})} \big)$, $\{ \mathbf{X}_{t},\mathbf{Y}_{t} \}_{t\in \llbracket 0 ; N_{0}\rrbracket}$ is
an $\{\mathbfcal{F}_{t}\}_{t\in \llbracket 0 ; N_{0}\rrbracket}$-adapted Markov chain with transition kernel $R$ and with initial distribution $\delta_{(\mathbf{x},\mathbf{y})}$.
Introduce the $ \llbracket 0 ; N_{0}\rrbracket$-valued random variable $\tau$ defined by
$$\tau=\begin{cases}
\inf\{k\in \llbracket 0 ; N_{0}-1\rrbracket : \mathbf{A}_{k}=1 \} & \text{if } \{k\in \llbracket 0 ; N_{0}-1\rrbracket : \mathbf{X}_{k}=1 \}\neq \emptyset, \\
N_{0} & \text{otherwise.}
\end{cases}
$$
It follows from Lemma \ref{algebra=} that $\tau$ is a stopping time with respect to
$\big\{\sigma\{\mathbf{Y}_{0},\ldots,\mathbf{Y}_{t}\}_{t\in \llbracket 0 ; N_{0}\rrbracket}\big\}$
showing that the control $\lambda$ defined by $\Big(\mathbf{\Omega},\mathbfcal{F},\mathbf{P}^{\pi}_{(\mathbf{x},\mathbf{y})}, 
\{\mathbfcal{F}_{t}\}_{t\in \llbracket 0 ; N_{0}\rrbracket},\{ \mathbf{X}_{t},\mathbf{Y}_{t} \}_{t\in \llbracket 0 ; N_{0}\rrbracket},\tau\Big)$ belongs to $\Lambda$.
Recalling that $\mathbf{Z}_{0}=0$ and that $\mathbf{Z}_{t}=I_{\{\mathbf{A}_{t-1}=1\}}+\mathbf{Z}_{t-1}I_{\{\mathbf{A}_{t-1}=0\}}$ for any $t\in \llbracket 1 ; N_{0} \rrbracket$, we get
that $\{\tau=t\}=\{\mathbf{Z}_{t}=0\}\cup\{\mathbf{A}_{t}=1\}$ for $t\in \llbracket 0 ; N_{0}-1 \rrbracket$ and $\{\tau=N_{0}\}=\{\mathbf{Z}_{N_{0}}=0\}$.
Now, by using the definitions of $H$ and $h$ it follows that
\begin{align*}
\sum_{t=0}^{N_{0}-1}  \mathbf{E}^{\pi}_{(\mathbf{x},\mathbf{y})} \big[ & H(\mathbf{X}_{t},\mathbf{Y}_{t},\mathbf{Z}_{t},\mathbf{A}_{t}) \big]
+\mathbf{E}^{\pi}_{(\mathbf{x},\mathbf{y})}\big[ h(\mathbf{X}_{N_{0}},\mathbf{Y}_{N_{0}},\mathbf{Z}_{N_{0}}) \big]
\nonumber \\
& = \sum_{t=0}^{N_{0}-1}  \mathbf{E}^{\pi}_{(\mathbf{x},\mathbf{y})} \big[ \mathbf{H}(\mathbf{X}_{t},\mathbf{Y}_{t}) I_{\{(\mathbf{Z}_{t},\mathbf{A}_{t})=(0,1)\}}\big]
+ \mathbf{E}^{\pi}_{(\mathbf{x},\mathbf{y})} \big[ \mathbf{H}(\mathbf{X}_{N_{0}},\mathbf{Y}_{N_{0}}) I_{\{\mathbf{Z}_{N_{0}}=1\}}\big]
\nonumber \\
& = \mathbf{E}^{\pi}_{(\mathbf{x},\mathbf{y})} \big[ \mathbf{H}(\mathbf{X}_{\tau},\mathbf{Y}_{\tau}) \big],
\end{align*}
implying that $\mathbfcal{H}(\mathbf{x},\mathbf{y},\ell)=\mathbfcal{H}_{\mathbfcal{M}}(\mathbf{x},\mathbf{y},\pi)$ and showing the second claim.
\hfill $\Box$

\bigskip

\noindent
\textbf{Proof of Theorem \ref{Equivalence-Costs}}
From Theorem 5.3.2 in \cite{bauerle11} we get that
$\overline{\mathbfcal{H}}_{\mathbfcal{M}}(\mathbf{x},\mathbf{y})=\overline{\mathbfcal{H}}_{\mathcal{M}}(\mathbf{x},\mathbf{y})$
and so from Proposition \ref{Equiv-partially-complete}, it follows that $\overline{\mathbfcal{H}}(\mathbf{x},\mathbf{y})=\overline{\mathbfcal{H}}_{\mathcal{M}}(\mathbf{x},\mathbf{y})$
giving the first equality in equation (\ref{Cost=}).
Under Assumptions \ref{HA1}, \ref{HB} and \ref{HD}, the hypotheses of Theorems 5.3.3 in \cite{bauerle11} are satisfied.
Therefore, it follows that
the Bellman equation $\{v_{k}\}_{ k \in \llbracket 0 ; N_{0}\rrbracket}$ for the model $\mathcal{M}$ is given by
\begin{align*}
\begin{cases}
v_{0}(\theta,y,z) = h(\theta,y,z) \\
v_{k}(\theta,y,z) = \max_{a\in \mathbb{A}} \big\{ H(\theta,y,z,a)+Qv_{k-1}(\theta,y,z,a)  \big\}
\end{cases}
\end{align*}
and satisfies $\overline{\mathbfcal{H}}_{\mathcal{M}}(\mathbf{x},\mathbf{y})=v_{N_{0}}(\delta_{\mathbf{x}},\mathbf{y},0)$.
However, since $h(\theta,y,1)=H(\theta,y,1)=0$, it easy to show that $v_{k}(\theta,y,1)=0$ for any
$(\theta,y)\in \mathcal{P}(\mathbb{X})\times\mathbb{Y}$ and $k\in \llbracket 0 ; N_{0}\rrbracket$.
Moreover, by using the definitions of $h$, $H$ and the kernel $Q$ we obtain
that $v_{0}(\theta,y,0) = \mathbf{H}(\theta,y)$ and 
\begin{align*}
v_{k}(\theta,y,0) & = \max_{a\in \mathbb{A}} \big\{ H(\theta,y,0,a)+Qv_{k-1}(\theta,y,0,a)  \big\} \nonumber \\
& = \max \big\{ \mathbf{H}(\theta,y),Sv_{k-1}(\theta,y)  \big\} = \mathfrak{B}v_{k-1}(\theta,y)
\end{align*}
for any $(\theta,y)\in \mathcal{P}(\mathbb{X})\times\mathbb{Y}$ and $k\in \llbracket 1 ; N_{0}\rrbracket$ implying that 
$\overline{\mathbfcal{H}}_{\mathcal{M}}(\mathbf{x},\mathbf{y})=\mathfrak{B}^{N_{0}}\mathbf{H}(\delta_{\mathbf{x}},\mathbf{y})$ and giving the second equality in equation (\ref{Cost=}).
\hfill $\Box$

\end{document}